\documentclass[english]{article}
\usepackage[T1]{fontenc}
\usepackage[latin9]{inputenc}
\usepackage{amsmath}
\usepackage{amssymb}
\usepackage{esint}
\usepackage{babel}
\begin{document}
\textbf{General monotonicity, interpolation of operators, and applications.}

{\begin{center} S.M. Grigoriev, Y. Sagher, T.R. Savage \end{center}}

\textbf{1. Introduction.}

$ $

\textbf{Definition 1.1}. A mesurable function $f$ on $\left(\Omega,\Sigma,\mu\right)$
is called \emph{admissible}, if the following conditions hold:

(1) $\left\{ f\neq0\right\} $ is a $\sigma$-finite set.

(2) $f<\infty$ a.e. on $\Omega$. 

(3) $\mu\left(\left\{ x:\left|f(x)\right|>\gamma\right\} \right)<\infty$
for some $\gamma>0$.

Recall the definition of the decreasing rearrangement of a function.

$ $

\textbf{Definition 1.2}. Let $(\Omega,\Sigma,\mu)$ be a measure space
and assume that $f$ is an admissible function on $\Omega$. Define,
for each $\alpha\geq0$, the \emph{distribution function} of $f$:
$f_{*}(\alpha)=\mu\left(\left\{ x:\left|f(x)\right|>\alpha\right\} \right)$. 

$ $

\textbf{Definition 1.3}. For each $x\geq0$, define 
\begin{equation}
f^{*}(x)=\begin{cases}
\begin{array}{cc}
\inf\left(\left\{ \alpha:f_{*}(\alpha)\leq x\right\} \right) & \qquad\mbox{if}\;\mbox{ \ensuremath{\exists}\ensuremath{\ensuremath{\alpha}:\ensuremath{f_{*}}(\ensuremath{\alpha})\ensuremath{\leq}x}}\\
\infty & \mbox{otherwise}
\end{array} & .\end{cases}\label{eq:f* - definition}
\end{equation}
$f^{*}$ is called the\emph{ decreasing rearrangement} of $f$.

The distribution function of $f^{*}$ with respect to Lebesgue measure
on $\mathbb{R}^{+}$ is equal to the distribution function of $f$
on $\left(\Omega,\Sigma,\mu\right)$:

\textbf{Theorem 1.4}. $\left(f^{*}\right)_{*}(\alpha)=f_{*}(\alpha),\forall\alpha\geq0$.

$ $

\textbf{Definition 1.5}. Considering sequences as functions defined
on $\Omega\subseteq\mathbb{Z}$, with $\mu=\#$, the counting measure,
we use formula $(\ref{eq:f* - definition})$ to define for each $k\in\mathbb{N}$,
\begin{equation}
a_{k}^{*}=\left(\left\{ a_{n}\right\} ^{*}\right)(k-1).\label{eq:Definition of a*}
\end{equation}

We recall the definition of weighted $L^{q}$-norms.

\textbf{Definition 1.6}. Let $0<q\leq\infty$, and let $w$ be a nonnegative
measurable function on $(\Omega,\Sigma,\mu)$. Assume that $f$ is
a measurable function on $(\Omega,\Sigma,\mu)$. Define the \emph{weighted
$L^{q}$-norm} of \emph{$f$} with the weight $w$:
\[
\left\Vert f\right\Vert _{L_{w}^{q}(\Omega,\Sigma,\mu)}=\left\Vert f\cdot w\right\Vert _{L^{q}(\Omega,\Sigma,\mu)}.
\]

$L_{w}^{q}=L_{w}^{q}(\Omega,\Sigma,\mu)$ is the space of all measurable
functions $f$ for which $\left\Vert f\right\Vert _{L_{w}^{q}(\Omega,\Sigma,\mu)}$
is finite.

In particular, working with $(\mathbb{N},2^{\mathbb{N}},\#)$, we
follow the practice of writing
\[
\left\Vert \left\{ a_{n}\right\} \right\Vert _{l_{w}^{q}}=\left\Vert \left\{ a_{n}\cdot w(n)\right\} \right\Vert _{l^{q}}.
\]
$l_{w}^{q}$ is the space of sequences $\left\{ a_{n}\right\} $ for
which $\left\Vert \left\{ a_{n}\right\} \right\Vert _{l_{w}^{q}}$
is finite.

We will make extensive use of the spaces of functions defined on $\mathbb{R}^{+}$
with the following weight:

\textbf{Definition 1.7}. Let $0<p\leq\infty$, $0<q\leq\infty$. We
denote, $\forall x>0$:
\begin{equation}
w(p,q)(x)=x^{\frac{1}{p}-\frac{1}{q}}.\label{eq: w(p,q) - definition}
\end{equation}

When we work with sequences, we denote, $\forall n\geq1$:
\begin{equation}
w(p,q)(n)=n^{\frac{1}{p}-\frac{1}{q}}.\label{eq: seq. w(p,q) - definition}
\end{equation}

\textbf{Definition 1.8}.

(i) Assume that $f$ is a measurable function on $(\Omega,\Sigma,\mu)$
and that $f^{*}<\infty$ on $\mathbb{R}^{+}$. For $0<p<\infty$,
$0<q\leq\infty$, and $p=q=\infty$, define the $L(p,q)$\emph{-norm
}of \emph{$f$}: 
\begin{equation}
\left\Vert f\right\Vert _{L(p,q)}=\left\Vert f\right\Vert _{L(p,q)(\Omega,\Sigma,\mu)}=\left\Vert f^{*}\right\Vert _{L_{w(p,q)}^{q}\left(\mathbb{R}^{+},{\cal B},\lambda\right)}\label{eq:L(p,q) - definition}
\end{equation}
$L(p,q)=L(p,q)(\Omega,\Sigma,\mu)$ is the space of all measurable
functions $f$, for which $\left\Vert f\right\Vert _{L(p,q)}$ is
finite.

(ii) Analogously,
\begin{equation}
\left\Vert \left\{ a_{n}\right\} \right\Vert _{l(p,q)}=\left\Vert \left\{ a_{n}^{*}\right\} \right\Vert _{l_{w(p,q)}^{q}}.\label{eq:l(p,q) - definition}
\end{equation}
$l(p,q)$ is the collection of all sequences $\left\{ a_{n}\right\} $
for which $\left\Vert \left\{ a_{n}\right\} \right\Vert _{l(p,q)}$
is finite.

$ $

\textbf{Lemma 1.9}. Assume that $\left\{ a_{n}\right\} ^{*}$ exists.
For $0<p,q\leq\infty$, 
\begin{equation}
\left\Vert \left\{ a_{k}\right\} \right\Vert _{l\left(p,q\right)}\sim\left\Vert \left\{ a_{n}\right\} ^{*}\left(\cdot\right)\right\Vert _{L_{w\left(p,q\right)}^{q}}.\label{eq:lpq and Lpq are similar for a and its function}
\end{equation}

$ $

\textbf{Definition 1.10}. $ONB$ is the collection of is an orthonormal
uniformly bounded sequences $\left\{ \phi_{n}\right\} $ of complex-valued
measurable functions.

$ $

\textbf{Theorem 1.11} (R.E.A.C. Paley, {[}10{]}).

Let $\left\{ \phi_{n}\right\} $ be an $ONB$ on $(\Omega,\Sigma,\mu)$,
$1<p\leq2$, and $c_{n}=\intop_{\Omega}f\bar{\phi_{n}}d\mu$:

(i) If $f\in L^{p}(\Omega,\Sigma,\mu)$, then $\left\Vert \left\{ c_{n}\right\} \right\Vert _{l(p',p)}\leq B\left\Vert f\right\Vert _{L^{p}}$.

(ii) If $\left\Vert \left\{ c_{n}\right\} \right\Vert _{l(p,p')}<\infty$
then $\left\Vert f\right\Vert _{L^{p'}}\leq B'\left\Vert \left\{ c_{n}\right\} \right\Vert _{l(p,p')}$.

The following classical theorem shows that, in the context of $L(p,q)$
spaces, Paley's theorem is best possible.

$ $

\textbf{Theorem 1.12} (G.H. Hardy and J.E. Littlewood, {[}6{]}).

Let $\left\{ c_{n}\right\} \searrow0$, $f(x)={\displaystyle \sum_{n=0}^{\infty}}c_{n}\cos nx$
or $f(x)={\displaystyle \sum_{n=1}^{\infty}}c_{n}\sin nx$, $p>1$.
Then $\left\Vert f\right\Vert _{L^{p}(0,\pi)}\sim\left\Vert \left\{ c_{n}\right\} \right\Vert _{l(p',p)}.$

Hardy and Littlewood also described a norm relation between monotone
decreasing functions and their trigonometric Fourier coefficients:

$ $

\textbf{Theorem 1.13} (G.H. Hardy and J.E. Littlewood, {[}7{]}).

Let $f\searrow$ on $[0,\pi]$, $f\geq0$, $1<p<\infty$, and let
$\left\{ c_{n}\right\} $ be the sequence of trigonometric Fourier
coefficients of $f$. Then: 
\[
\left\Vert f\right\Vert _{L^{p}(0,\pi)}\sim\left|c_{0}\right|+\left\Vert \left\{ c_{n}\right\} \right\Vert _{l_{w(p',p)}^{p}.}
\]

The extension of Theorem 1.13 to broader classes of sequences started
with the introduction of \emph{quasi-monotone} \emph{sequences} by
Shah {[}13{]}; this concept also applies to functions defined on $(0,\infty)$:

$ $

\textbf{Definition 1.14}. Let $\beta>0$.

(a) If $a_{n}\geq0$ and $\left\{ a_{n}\cdot n^{-\beta}\right\} \searrow$,
the sequence $\left\{ a_{n}\right\} $ is called a \emph{$\beta$-quasi-monotone
decreasing sequence} (written as $\left\{ a_{n}\right\} \in QDS\left(\beta\right)$).

(b) If $f:\mathbb{R}^{+}\rightarrow\mathbb{R}^{+}$ is such that $x^{-\beta}f(x)\searrow$
then $f$ is called a \emph{$\beta$-quasi-monotone decreasing function},
$f\in QD\left(\beta\right)$.

Denote $QDS=\underset{\beta>0}{\cup}QDS\left(\beta\right)$; $QD=\underset{\beta>0}{\cup}QD\left(\beta\right)$.
We say that $\left\{ a_{n}\right\} $ is a \emph{quasi-monotone decreasing
sequence} if $\left\{ a_{n}\right\} \in QDS$ and that $f$ is a \emph{quasi-monotone
decreasing function} if $f\in QD$.

The interest in quasi-monotone sequences stems from the fact that
the class of trigonometric series with quasi-monotone Fourier coefficients
is closed under term-by-term differentiation. For an application,
see {[}11{]}.

Askey and Wainger {[}1{]} proved that if $\left\{ c_{n}\right\} $
is a quasi-monotone sequence of trigonometric Fourier coefficients
of $f$ then for $1<p<\infty$, $1\leq q<\infty,$ 
\begin{equation}
\left\Vert f\right\Vert _{L_{w(p,q)(0,\pi)}^{q}}\sim\left\Vert \left\{ c_{n}\right\} \right\Vert _{l_{w(p',q)}^{q}}.\label{eq:Askey-Wainger}
\end{equation}

An extension of Real Interpolation theory made it possible to simplify
the proof of $(\ref{eq:Askey-Wainger})$. The extension consists of
interpolation of normed \emph{monoids} rather than normed vector spaces.

Tikhonov {[}14{]} defined the following extension of the class of
quasi-monotone sequences:

\textbf{Definition} \textbf{1.15}. Let $B\geq1$. We write $\left\{ a_{n}\right\} \in GMS\left(B\right)$
if for all $n\geq1$: 
\begin{equation}
\sum_{k=n}^{2n-1}\left|a_{k}-a_{k+1}\right|\leq B\left|a_{n}\right|.\label{eq:Definition of GMS}
\end{equation}

Define $GMS=\underset{B\geq1}{\cup}GMS(B)$. If $\left\{ a_{n}\right\} \in GMS$,
we say that $\left\{ a_{n}\right\} $ is a \emph{general monotone}
\emph{sequence. }

Liflyand and Tikhonov {[}9{]} defined \emph{general monotone functions}
($GM$). We will work with a special case of their definition.

\textbf{Definition 1.16}. Let \textbf{$B\geq1$. }We write $f\in GM\left(B\right)$
if $f\in BV_{loc}\left(\left(0,\infty\right)\right)$ and for all
$x>0$: 
\begin{equation}
V_{f}\left(\left[x,2x\right]\right)\leq B\left|f(x)\right|,\label{eq:GMF: definition}
\end{equation}
where $V_{f}\left(\left[a,b\right]\right)$ is the variation of $f$
on $\left[a,b\right]$. Define $GM=\underset{B\geq1}{\cup}GM(B)$.

If $f\in GM$, we say that $f$ is a \emph{general monotone function.} 

We will often use the following classes of functions and sequences:

\textbf{Definition 1.17}. 
\begin{equation}
GM_{1}\left(B\right)=\left\{ f\in{\cal M}_{\mathbb{R}^{+}}:\left|f(t)\right|\leq B\left|f(x)\right|,\mbox{ }\forall x\leq t\leq2x\right\} ;\label{eq:almost monotone-1}
\end{equation}
\begin{equation}
GMS_{1}\left(B\right)=\left\{ \left\{ a_{k}\right\} _{k=1}^{\infty}:\left|a_{k}\right|\leq B\left|a_{n}\right|,\mbox{ }\forall n\leq k\leq2n\right\} ;\label{eq:almost monotone}
\end{equation}
\begin{equation}
{\scriptscriptstyle {\scriptstyle {\textstyle GM_{2}\left(B\right)=}}\left\{ {\scriptstyle f\in BV_{loc}\left(\mathbb{R}^{+}\right):V_{f}([x,M])\leq B\left(|f(x)|+{\displaystyle \intop_{x}^{M}}\left|f(t)\right|{\displaystyle \frac{dt}{t}}\right),\;\forall M>0,\forall0<x<M}\right\} ;}\label{eq:almost regular GM-1}
\end{equation}
\begin{equation}
{\textstyle {\textstyle GMS_{2}\left(B\right)}=\left\{ {\scriptstyle \left\{ a_{k}\right\} _{k=1}^{\infty}:{\displaystyle \sum_{k=n}^{N-1}}\left|a_{k}-a_{k+1}\right|\leq B\left(|a_{n}|+{\displaystyle \sum_{k=n+1}^{N}}{\textstyle \frac{|a_{k}|}{k}}\right)},{\scriptstyle \forall N>1,\forall1\leq n<N}\right\} }.\label{eq:almost regular GM}
\end{equation}
We write for $j=1,2$, $GM_{j}=\underset{B\geq1}{\cup}GM_{j}\left(B\right)$,
$GM_{j}^{+}=\left\{ f\in GM_{j}:f\geq0\right\} $, $GMS_{j}=\underset{B\geq1}{\cup}GMS_{j}\left(B\right)$,
and $GMS_{j}^{+}=GMS_{j}\cap\left\{ \left\{ a_{n}\right\} :a_{n}\geq0,\forall n\right\} $.

$ $

\textbf{Definition 1.18}. For each $0\leq\varphi<{\displaystyle \frac{\pi}{2}}$
and $\alpha\in\mathbb{R}$, define 
\begin{equation}
S_{\alpha,\varphi}:=\left\{ z\in\mathbb{C}:\left|\arg\left(e^{-i\alpha}z\right)\right|\leq\varphi\right\} \cup\left\{ 0\right\} .\label{eq:cone in C}
\end{equation}

\textbf{Lemma 1.19}. Let $(\Omega,\Sigma,\mu)$ be a measure space
and assume that $f(\omega)\in S_{\alpha,\varphi}$ for a.e. $\omega\in\Omega$.
Then
\begin{equation}
\intop_{\Omega}\left|f\right|d\mu\leq\frac{1}{\cos\varphi}\left|\intop_{\Omega}fd\mu\right|.\label{eq:reverse triang. ineq. for integrals in S}
\end{equation}

\textbf{Definition 1.20}. For each $0\leq\varphi<{\displaystyle \frac{\pi}{2}}$
and $0\leq\alpha<2\pi$, define
\[
GM_{\alpha,\varphi}=\left\{ f\in GM:f(x)\in S_{\alpha,\varphi},\mbox{ }\mu\mbox{-a.e. }x>0\right\} ;
\]
\[
GMS_{\alpha,\varphi}=\left\{ \left\{ a_{n}\right\} \in GMS:a_{n}\in S_{\alpha,\varphi},\forall n\geq1\right\} .
\]
We write $GM_{\alpha,\varphi}(B)=GM_{\alpha,\varphi}\cap GM(B)$ and
$GMS_{\alpha,\varphi}(B)=GMS_{\alpha,\varphi}\cap GMS(B)$.

$ $

\textbf{Remark 1.21}. Observe that $GMS_{0,0}=\left\{ \left\{ a_{n}\right\} \in GMS:a_{n}\geq0\right\} =:GMS^{+}$,
$GM_{0,0}=\left\{ f\in GM:f\geq0\right\} =:GM^{+}$.

$ $

\textbf{Lemma 1.22} (E. Liflyand and S. Tikhonov, {[}9{]}).

\begin{equation}
GM\left(B\right)\subset GM_{1}\left(2B\right)\cap GM_{2}\left(2B^{2}\right);\label{eq:almost monotone-1: B}
\end{equation}
\begin{equation}
GM_{1}\left(B\right)\cap GM_{2}\left(B\right)\subset GM\left(2B^{2}\right).\label{eq:almost regular-1: B}
\end{equation}

Analogously for general monotone sequences:

\textbf{Lemma 1.23} (S. Tikhonov, {[}14{]}).
\begin{equation}
GMS(B)\subset GMS_{1}(2B)\cap GMS_{2}\left(2B^{2}\right);\label{eq:almost monotone: B}
\end{equation}
\begin{equation}
GMS_{1}(B)\cap GMS_{2}(B)\subset GMS\left(2B^{2}\right).\label{eq:almost regular GM: B}
\end{equation}

\newpage

\textbf{2. Some results from Real Interpolation theory.}

$ $

\textbf{Definition} \textbf{2.1}. A \emph{semigroup }in a topological
vector space (t.v.s.)\emph{ }is a subset $G$ of the t.v.s. that is
closed under addition. A semigroup that contains $0$ is a \emph{monoid}.
A monoid that is closed under multiplication by any positive scalar
is a \emph{cone}.

\textbf{Definition} \textbf{2.2}. A \emph{quasi-norm} on a vector
space $X$ is a function $\left\Vert \cdot\right\Vert _{X}:X\rightarrow[0,\infty)$
that satisfies:

(a) $\left\Vert a\right\Vert _{X}=0\Leftrightarrow a=0$.

(b) $\exists k=k(X)$ such that $\forall a_{1},a_{2}\in X$: $\left\Vert a_{1}+a_{2}\right\Vert _{X}\leq k\left(\left\Vert a_{1}\right\Vert _{X}+\left\Vert a_{2}\right\Vert _{X}\right).$

(c) $\left\Vert \lambda a\right\Vert _{G}=\left|\lambda\right|\left\Vert a\right\Vert _{G}$,
$\forall\lambda\in\mathbb{C},a\in G$.

A vector space equipped with a quasi-norm is a \emph{quasi-normed
vector space}, a complete quasi-normed vector space is called a \emph{quasi-Banach
space}.

$ $

\textbf{Definition 2.3}. If $G_{0}$ and $G_{1}$ are monoids continuously
embedded in a t.v.s., we say that $(G_{0},G_{1})$ is an \emph{interpolation
couple}.

Throughout this section, $G_{0}$, $G_{1}$ will be quasi-normed monoids
which form an interpolation couple in a t.v.s. of measurable functions,
the topology of the t.v.s. is defined by the convergence in measure.

\textbf{Definition} \textbf{2.4}. For $f\in G_{0}+G_{1}$, $t>0$,
define 
\[
K(t,f)=K(t,f,G_{0},G_{1})
\]
\begin{equation}
=\inf\left\{ \left.\left\Vert f_{0}\right\Vert _{G_{0}}+t\left\Vert f_{1}\right\Vert _{G_{1}}\right|f_{0}+f_{1}=f,f_{j}\in G_{j},j=0,1\right\} .\label{eq:K-functional}
\end{equation}

$ $

\textbf{Definition} \textbf{2.5}. ($K$-method of interpolation of
monoids).

For $f\in G_{0}+G_{1}$, and for each $0<\theta<1$ and $0<q\leq\infty$
or for $0\leq\theta\leq1$ and $q=\infty$, define
\begin{equation}
\left\Vert f\right\Vert _{\left(G_{0},G_{1}\right)_{\theta,q;K}}=\left\Vert t^{-\theta}K(t,f)\right\Vert _{L^{q}\left((0,\infty),\frac{dt}{t}\right)}.\label{eq: K-method}
\end{equation}

Define the \emph{Interpolation monoid with respect to $K$} as
\begin{equation}
(G_{0},G_{1})_{\theta,q;K}=\left\{ f\in G_{0}+G_{1}\left|\left\Vert f\right\Vert _{\left(G_{0},G_{1}\right)_{\theta,q;K}}<\infty\right.\right\} .\label{eq:K-interpolation space}
\end{equation}

It is clear that $(G_{0},G_{1})_{\theta,q;K}$ is a monoid. When $G_{0}$
and $G_{1}$ are quasi-normed vector spaces, we call $(G_{0},G_{1})_{\theta,q;K}$
an \emph{Interpolation space}.

$ $

\textbf{Theorem 2.6} (J. Gilbert, {[}5{]}).

Assume $X$ is a quasi-normed monoid, $w$ is a measurable positive
function on $\Omega$, $\sigma$ is positive, piecewise continuous
function such that both $\sigma$ and $\sigma_{t}$, where $\sigma_{t}(\lambda):=t\lambda\sigma(t\lambda)$,
are in $L^{\infty}(0,\infty)$, $\forall t>0$. For any $f\in(X,X_{w})_{\theta,q;K}$,
$\theta\in(0,1)$, and $q\in(0,\infty]$, 
\begin{equation}
\left\Vert f\right\Vert _{\theta,q;K}\sim\left(\intop_{0}^{\infty}\left(t^{-\theta}\left\Vert f\cdot\sigma_{t}\circ w\right\Vert _{X}\right)^{q}\frac{dt}{t}\right)^{1/q}.\label{eq:Gilbert's identity}
\end{equation}

\newpage

\textbf{3. Selected results from the theory of $L(p,q)$ spaces.}

$ $

\textbf{Theorem 3.1} ({[}8{]}). Suppose that $I$ is a subinterval
of $\mathbb{Z}$ and that $\left\{ a_{k}\right\} _{k\in I}$, $\left\{ c_{k}\right\} _{k\in I}$
are nonnegative sequences. If $\left\{ c_{k}\right\} \searrow$ then
${\displaystyle \sum_{k\in I}a_{k}c_{k}\leq}{\displaystyle \sum_{k\in I}}a_{k}^{*}c_{k}.$
If $\left\{ c_{k}\right\} \nearrow$ then ${\displaystyle \sum_{k\in I}a_{k}c_{k}\geq}{\displaystyle \sum_{k\in I}}a_{k}^{*}c_{k}.$

$ $

\textbf{Theorem 3.2} (G.H. Hardy).

Let $f_{1},f_{2},g$ be nonnegative, measurable functions on $\mathbb{R}^{+}$.

(a) Assume that $g\searrow$ and that for all $a>0$, ${\displaystyle \intop_{0}^{a}}f_{1}dt\leq{\displaystyle \intop_{0}^{a}}f_{2}dt$.
Then
\begin{equation}
{\displaystyle \intop_{0}^{\infty}}f_{1}gdt\leq{\displaystyle \intop_{0}^{\infty}}f_{2}gdt.\label{eq:rearrangement - general - integrall}
\end{equation}
(b) Assume that $g\nearrow$ and that for all $a>0$, ${\displaystyle \intop_{a}^{\infty}}f_{1}dt\leq{\displaystyle \intop_{a}^{\infty}}f_{2}dt$.
Then $(\ref{eq:rearrangement - general - integrall})$ holds.

$ $

\textbf{Corollary} \textbf{3.3}. Suppose that $f,g$ are nonnegative,
measurable on $\mathbb{R}^{+}$ functions.

If $g\searrow$ then ${\displaystyle {\displaystyle \intop_{0}^{\infty}}fgdt\leq}{\displaystyle \intop_{0}^{\infty}}f^{*}gdt.$
If $g\nearrow$ then ${\displaystyle {\displaystyle \intop_{0}^{\infty}}fgdt\geq}{\displaystyle \intop_{0}^{\infty}}f^{*}gdt.$

$ $

\textbf{Theorem 3.4}. Let $f$ be measurable on $\mathbb{R}^{+}$.
For $0<p\leq q\leq\infty$, 
\begin{equation}
\left\Vert f\right\Vert _{L(p,q)\left(0,\infty\right)}\leq\left\Vert f\right\Vert _{L_{w(p,q)\left(0,\infty\right)}^{q}}.\label{eq: a simple case for Lpqs}
\end{equation}

$ $

The following theorem is proved for sine Fourier coefficients in {[}12{]}.

$ $

\textbf{Theorem 3.5.} Let $f\in L(p',q)$, $1<p<\infty$, $0<q\leq\infty$
and assume that $\left\{ c_{n}\right\} $ are Fourier coefficients
of $f$ with respect to $\left\{ e^{int}\right\} _{n=-\infty}^{\infty}$.
Define $\sigma_{n}={\displaystyle \frac{1}{2n+1}}{\displaystyle \sum_{k=-n}^{n}}c_{k}$,
and let $\left\{ m\left(a_{n}\right)\right\} $ be the least bell-shaped
majorant of $\left\{ a_{n}\right\} $, that is to say, $m\left(a_{n}\right):=\sup\left\{ \left|a_{k}\right|:\left|k\right|\geq\left|n\right|\right\} $.
Then $\exists A(p)>0$ so that: 
\begin{equation}
\left\Vert \left\{ m\left(\sigma_{n}\right)\right\} \right\Vert _{l(p,q)}\leq A(p)\left\Vert f\right\Vert _{L(p',q)(0,2\pi)},\label{eq:Sagher-12}
\end{equation}

Consequently,

\begin{equation}
\left\Vert \left\{ \sigma_{n}\right\} \right\Vert _{l(p,q)}\leq A(p)\left\Vert f\right\Vert _{L(p',q)(0,2\pi)}.\label{eq: Sagher-12-1}
\end{equation}

\newpage

\textbf{4.} \textbf{Some elementary properties of $GM$ and $GMS$. }

$ $

\textbf{Lemma 4.1}. For $j=1,2$, if $f\in GM_{j}$ then $\left\{ f(n)\right\} \in GMS_{j}$.
Moreover, if $\left\{ a_{n}\right\} \in GMS_{1}$ and $f(x)=a_{\left\lceil x\right\rceil }$
then $f\in GM_{1}$.

\textbf{Lemma 4.2}. If $f\in GM$ then $\left\{ f(n)\right\} \in GMS$.
Conversely, for each $\left\{ a_{n}\right\} \in GMS$ there is an
$f\in GM$ such that $f(n)=a_{n}$, $\forall n\geq1$.

Using Lemma 4.2, results for $GMS$ follow from corresponding results
for $GM$, and therefore we shall only give proofs for $GM$ functions.

$ $

\textbf{Lemma 4.3}. Assume that $f\in GM_{1}\left(B\right)$ and that
$f^{*}$ exists. Then $\forall x>0$,
\begin{equation}
\left|f(x)\right|\leq Bf^{*}\left(\left(\frac{x}{2}\right)-\right).\label{eq:f and f* for GM}
\end{equation}

\textbf{Proof}.

Let $x>0$, $0<\theta<1$. Since $(\ref{eq:f and f* for GM})$ trivially
holds if $f(x)=0$, we can assume that $f(x)\neq0$. For all $t\in\left[\frac{x}{2},x\right]$,
$t\leq x\leq2t$, and by $(\ref{eq:almost monotone-1})$,
\[
\left|f(t)\right|\geq{\displaystyle \frac{\left|f(x)\right|}{B}}>{\displaystyle \frac{\theta\left|f(x)\right|}{B}}.
\]
Therefore:
\[
f_{*}\left(\frac{\theta\left|f(x)\right|}{B}\right)=\lambda\left\{ t:\left|f(t)\right|>{\displaystyle \frac{\theta\left|f(x)\right|}{B}}\right\} \geq\frac{x}{2}>\frac{\theta x}{2},
\]
and since $f_{*}\searrow$:
\[
f^{*}\left(\frac{\theta x}{2}\right)=\inf\left\{ \alpha:f_{*}(\alpha)\leq\frac{\theta x}{2}\right\} \geq\frac{\theta\left|f(x)\right|}{B},
\]
or $|f(x)|\leq\frac{B}{\theta}f^{*}\left(\frac{\theta x}{2}\right)$.
Letting $\theta\rightarrow1^{-}$, $(\ref{eq:f and f* for GM})$ follows.$\square$

The corresponding result for $GMS_{1}$ was proved in {[}3{]}:

\textbf{Lemma 4.4} (B. Booton, {[}3{]}).

Let $\left\{ c_{k}\right\} \in GMS_{1}\left(B\right)$ be such that
$\left\{ c_{k}^{*}\right\} $ exists. Then $\forall n\geq1$, 
\[
\left|c_{n}\right|\leq Bc_{\left\lfloor \frac{n}{2}\right\rfloor +1}^{*}.
\]

\textbf{Lemma 4.5}. Assume that $f\in GM_{1}$, $A>0$, $\alpha\in\mathbb{R}$.

(i) If ${\displaystyle \intop_{A}^{\infty}}\left|f(x)\right|x^{\alpha}dx<\infty$
then ${\displaystyle \lim_{x\rightarrow\infty}}x^{\alpha+1}f(x)=0$.

(ii) If ${\displaystyle \intop_{0}^{A}}\left|f(x)\right|x^{\alpha}dx<\infty$
then ${\displaystyle \lim_{x\rightarrow0+}}x^{\alpha+1}f(x)=0$.

\textbf{Proof}.

Let $2^{k}\leq x\leq2^{k+1}$. If $\alpha<0$ then $2^{\alpha}\cdot2^{k\alpha}=2^{\left(k+1\right)\alpha}\leq x^{\alpha}\leq2^{k\alpha}$.
If $\alpha\geq0$ then $2^{k\alpha}\leq x^{\alpha}\leq2^{\left(k+1\right)\alpha}=2^{\alpha}\cdot2^{k\alpha}$.
That is to say, for any $\alpha\in\mathbb{R}$: 
\[
\min\left\{ 2^{\alpha},1\right\} \cdot2^{k\alpha}\leq x^{\alpha}\leq\max\left\{ 2^{\alpha},1\right\} \cdot2^{k\alpha}.
\]
Furthermore, by $(\ref{eq:almost monotone-1})$, ${\displaystyle \frac{1}{B}}\left|f\left(2^{k+1}\right)\right|\leq\left|f(x)\right|\leq B\left|f\left(2^{k}\right)\right|$.
Therefore:
\begin{equation}
\frac{1}{B}\min\left\{ 2^{\alpha},1\right\} 2^{k\alpha}\left|f\left(2^{k+1}\right)\right|\leq x^{\alpha}\left|f(x)\right|\leq B\max\left\{ 2^{\alpha},1\right\} 2^{k\alpha}\left|f\left(2^{k}\right)\right|.\label{eq: boundedness of GM by dyadic values}
\end{equation}
It follows from $(\ref{eq: boundedness of GM by dyadic values})$
that for all $-\infty\leq n<N\leq\infty$, 
\[
{\displaystyle \intop_{2^{n}}^{2^{N}}}\left|f(x)\right|x^{\alpha}dx=\sum_{k=n}^{N-1}\intop_{2^{k}}^{2^{k+1}}x^{\alpha}\left|f(x)\right|dx
\]
\[
\geq\frac{\min\left\{ 2^{\alpha},1\right\} }{B}\cdot\sum_{k=n}^{N-1}2^{k\alpha}\left|f\left(2^{k+1}\right)\right|\cdot2^{k}=\frac{\min\left\{ 2^{\alpha},1\right\} }{B\cdot2^{\alpha+1}}\sum_{k=n+1}^{N}2^{k\left(\alpha+1\right)}\left|f\left(2^{k}\right)\right|.
\]
Therefore, if ${\displaystyle \intop_{A}^{\infty}}\left|f(x)\right|x^{\alpha}dx<\infty$
then ${\displaystyle \sum_{k=\left\lceil \log_{2}A\right\rceil }^{\infty}}2^{k\left(\alpha+1\right)}\left|f\left(2^{k}\right)\right|<\infty$,
and so ${\displaystyle \lim_{k\rightarrow\infty}}2^{k\left(\alpha+1\right)}\left|f\left(2^{k}\right)\right|=0.$
If ${\displaystyle \intop_{0}^{A}}\left|f(x)\right|x^{\alpha}dx<\infty$
then ${\displaystyle \sum_{k=-\infty}^{\left\lfloor \log_{2}A\right\rfloor }}2^{k\left(\alpha+1\right)}\left|f\left(2^{k}\right)\right|<\infty$,
and so ${\displaystyle \lim_{k\rightarrow-\infty}}2^{k\left(\alpha+1\right)}\left|f\left(2^{k}\right)\right|=0.$
The claim of the Lemma now follows by substituting $\alpha+1$ for
$\alpha$ in the right-hand side of $(\ref{eq: boundedness of GM by dyadic values})$.
$\square$

$ $

\textbf{Corollary} \textbf{4.5}. If $0<p\leq\infty$, $0<q<\infty$
and $f\in L_{w(p,q)}^{q}\cap GM_{1}$ then ${\displaystyle \lim_{x\rightarrow0+}}x^{\frac{1}{p}}f(x)={\displaystyle \lim_{x\rightarrow\infty}}x^{\frac{1}{p}}f(x)=0$.

By Lemma 4.1:

\textbf{Lemma 4.6}. Assume $\left\{ a_{n}\right\} \in GMS_{1}$. If
$\alpha\in\mathbb{R}$ is such that ${\displaystyle \sum_{n=1}^{\infty}}n^{\alpha}\left|a_{n}\right|<\infty$
then ${\displaystyle \lim_{n\rightarrow\infty}}n^{\alpha+1}a_{n}=0$.

\textbf{Corollary} \textbf{4.6}. If $0<p\leq\infty$, $0<q<\infty$
and $\left\{ a_{n}\right\} \in l_{w(p,q)}^{q}\cap GMS_{1}$ then ${\displaystyle \lim_{n\rightarrow\infty}}n^{\frac{1}{p}}a_{n}=0$.

If $\alpha=-1$ then the claim of Lemma 4.6 is also true for $GMS_{2}$:

$ $

\textbf{Lemma 4.7}. Assume $\left\{ a_{n}\right\} \in l_{\frac{1}{k}}^{1}\cap GMS_{2}$.
Then ${\displaystyle \lim_{n\rightarrow\infty}}a_{n}=0$.

\textbf{Proof.}

By the assumption of the Lemma, ${\displaystyle \sum_{m=1}^{\infty}}\left|a_{m}-a_{m+1}\right|\leq B{\displaystyle \sum_{k=1}^{\infty}}{\displaystyle \frac{|a_{k}|}{k}}<\infty$,
and so, ${\displaystyle \lim_{m\rightarrow\infty}}\left(a_{1}-a_{m}\right)={\displaystyle \sum_{j=1}^{\infty}}\left(a_{m}-a_{m+1}\right)$
exists and finite. If $0\neq a={\displaystyle \lim_{m\rightarrow\infty}}a_{m}$
then $\exists N$ such that $\left|\left|a_{m}\right|-\left|a\right|\right|<{\displaystyle \frac{\left|a\right|}{2}}$,
$\forall m\geq N$, that is, ${\displaystyle \frac{\left|a\right|}{2m}}\leq{\displaystyle \frac{\left|a_{m}\right|}{m}}$.
But then ${\displaystyle \sum_{m=N}^{\infty}\frac{\left|a\right|}{2m}<\infty}$,
since $\left\{ a_{n}\right\} \in l_{\frac{1}{k}}^{1}$, a contradiction$.\square$

$ $

\textbf{Lemma 4.8}. Let $B\geq1$, $f\in GM_{1}\left(B\right)$. Then
for $0<p<\infty$, $0<q<\infty$: 
\begin{equation}
\frac{C_{1}(p,q)}{B}\left\Vert f\right\Vert _{L_{w\left(p,q\right)}^{q}}\leq\left(\sum_{k=-\infty}^{\infty}2^{\frac{kq}{p}}\left|f\left(2^{k}\right)\right|^{q}\right)^{\frac{1}{q}}\leq BC_{2}(p,q)\left\Vert f\right\Vert _{L\left(p,q\right)},\label{eq: another estimate between Lqw and}
\end{equation}
with $0<C_{j}(p,q)<\infty$. For $0<q<\infty$, 
\begin{equation}
\left\Vert f\right\Vert _{L_{w\left(\infty,q\right)}^{q}}=\left(\intop_{0}^{\infty}\left|f\left(x\right)\right|^{q}\frac{dx}{x}\right)^{\frac{1}{q}}\leq B\left(\ln2\right)^{\frac{1}{q}}\left(\sum_{k=-\infty}^{\infty}\left|f\left(2^{k}\right)\right|^{q}\right)^{\frac{1}{q}}.\label{eq:Lqw and the discrete norm: case of infinite p}
\end{equation}
Furthermore, for $0<p\leq\infty$:
\[
\frac{1}{2^{\frac{1}{p}}B}\left\Vert f\right\Vert _{L\left(p,\infty\right)}\leq\frac{1}{2^{\frac{1}{p}}B}\left\Vert f\right\Vert _{L_{w\left(p,\infty\right)}^{\infty}}\leq\sup_{k\in\mathbb{Z}}\left\{ 2^{\frac{k}{p}}\left|f\left(2^{k}\right)\right|\right\} 
\]
 
\begin{equation}
\leq2^{\frac{1}{p}}B\left\Vert f\right\Vert _{L_{w\left(p,\infty\right)}^{\infty}}\leq2^{\frac{2}{p}}B^{2}\left\Vert f\right\Vert _{L\left(p,\infty\right)}.\label{eq:estimate for GM for infinite q}
\end{equation}

\textbf{Proof}.

For $f\in GM_{1}\left(B\right)$ and $2^{k}\leq x\leq2^{k+1}$, $\left|f\left(x\right)\right|\leq B\left|f\left(2^{k}\right)\right|$,
and so:

\[
\left\Vert f\right\Vert _{L_{w\left(p,q\right)}^{q}}^{q}=\intop_{0}^{\infty}x^{\frac{q}{p}-1}\left|f\left(x\right)\right|^{q}dx=2^{\frac{q}{p}}\sum_{k=-\infty}^{\infty}\intop_{2^{k}}^{2^{k+1}}\left(\frac{x}{2}\right)^{\frac{q}{p}}\left|f\left(x\right)\right|^{q}\frac{dx}{x}
\]
\[
\leq2^{\frac{q}{p}}B^{q}\sum_{k=-\infty}^{\infty}2^{\frac{kq}{p}}\left|f\left(2^{k}\right)\right|^{q}\intop_{2^{k}}^{2^{k+1}}\frac{dx}{x}=2^{\frac{q}{p}}B^{q}\ln2{\displaystyle \sum_{k=-\infty}^{\infty}}2^{\frac{kq}{p}}\left|f\left(2^{k}\right)\right|^{q},
\]
proving the first inequality in $(\ref{eq: another estimate between Lqw and})$
and $(\ref{eq:Lqw and the discrete norm: case of infinite p})$. Assume
$2^{k-2}\leq x\leq2^{k-1}$, or equivalently, $2x\leq2^{k}\leq4x$.
If $q>p$ then $2^{k\left(\frac{q}{p}-1\right)}\leq\left(4x\right)^{\frac{q}{p}-1}$,
if $q\leq p$ then $2^{k\left(\frac{q}{p}-1\right)}\leq\left(2x\right)^{\frac{q}{p}-1}$
and so, for all $0<p,q<\infty$: 
\[
2^{k\left(\frac{q}{p}-1\right)}\leq A(p,q)x^{\frac{q}{p}-1},
\]
$A(p,q)=\max\left\{ 2^{\frac{q}{p}-1},4^{\frac{q}{p}-1}\right\} $.
By Lemma 4.3, $\left|f\left(2^{k}\right)\right|\leq Bf^{*}\left(2^{k-1}-\right)$,
and so:

\[
\sum_{k=-\infty}^{\infty}2^{\frac{kq}{p}}\left|f\left(2^{k}\right)\right|^{q}\leq B^{q}\sum_{k=-\infty}^{\infty}2^{\frac{kq}{p}}\left(f^{*}\left(2^{k-1}-\right)\right)^{q}
\]
\[
=4B^{q}\sum_{k=-\infty}^{\infty}\intop_{2^{k-2}}^{2^{k-1}}2^{k\left(\frac{q}{p}-1\right)}\left(f^{*}\left(2^{k-1}-\right)\right)^{q}dx\leq4A(p,q)B^{q}\intop_{0}^{\infty}x^{\frac{q}{p}-1}\left(f^{*}\left(x\right)\right)^{q}dx
\]
\[
=B^{q}A\left(p,q\right)\intop_{0}^{\infty}x^{\frac{q}{p}-1}\left(f^{*}\left(x\right)\right)^{q}dx=B^{q}A\left(p,q\right)\left\Vert f\right\Vert _{L\left(p,q\right)}^{q},
\]
proving the second inequality in $(\ref{eq: another estimate between Lqw and})$.

Let us consider $(\ref{eq:estimate for GM for infinite q})$. Observe
that $\left|f\left(x\right)\right|\leq x^{-\frac{1}{p}}{\displaystyle \sup_{x>0}}\left\{ x^{\frac{1}{p}}\left|f\left(x\right)\right|\right\} $,
and so, $\left|f^{*}\left(x\right)\right|\leq x^{-\frac{1}{p}}{\displaystyle \sup_{x>0}}\left\{ x^{\frac{1}{p}}\left|f\left(x\right)\right|\right\} $,
that is to say, $\left\Vert f\right\Vert _{L\left(p,\infty\right)}\leq\left\Vert f\right\Vert _{L_{w\left(p,\infty\right)}^{\infty}}$,
proving the first inequality. For $2^{k}\leq x\leq2^{k+1}$ and $\alpha={\displaystyle \frac{1}{p}}$,
$(\ref{eq: boundedness of GM by dyadic values})$ shows that: 
\[
\frac{1}{B}\cdot2^{\frac{k}{p}}\left|f\left(2^{k+1}\right)\right|\leq x^{\frac{1}{p}}\left|f(x)\right|\leq2^{\frac{1}{p}}B\cdot2^{\frac{k}{p}}\left|f\left(2^{k}\right)\right|,
\]
or equivalenty, using $2^{k-1}\leq{\displaystyle \frac{x}{2}}\leq2^{k}$:
\[
\frac{1}{2^{\frac{1}{p}}B}x^{\frac{1}{p}}\left|f(x)\right|\leq2^{\frac{k}{p}}\left|f\left(2^{k}\right)\right|\leq2^{\frac{1}{p}}B\cdot\left(\frac{x}{2}\right)^{\frac{1}{p}}\left|f\left(\frac{x}{2}\right)\right|.
\]
Taking the supremum over all $k\in\mathbb{Z}$ we obtain: 
\[
\frac{1}{2^{\frac{1}{p}}B}\left\Vert f\right\Vert _{L_{w\left(p,\infty\right)}^{\infty}}\leq\sup_{k\in\mathbb{Z}}\left\{ 2^{\frac{k}{p}}\left|f\left(2^{k}\right)\right|\right\} \leq2^{\frac{1}{p}}B\left\Vert f\right\Vert _{L_{w\left(p,\infty\right)}^{\infty}},
\]
proving the second and the third inequality. Finally, by Lemma 4.3,
for all $x>0$, $\left|f\left(x\right)\right|\leq Bf^{*}\left(\left(\frac{x}{2}\right)-\right)$,
therefore, for all $0<\theta<1$: 
\[
{\displaystyle \sup_{x>0}}\left\{ x^{\frac{1}{p}}\left|f\left(x\right)\right|\right\} \leq B\sup\left\{ x^{\frac{1}{p}}f^{*}\left(\left(\frac{x}{2}\right)-\right)\right\} =2^{\frac{1}{p}}B\sup\left\{ x^{\frac{1}{p}}f^{*}\left(x-\right)\right\} 
\]
\[
\leq2^{\frac{1}{p}}B{\displaystyle \sup_{x>0}}\left\{ x^{\frac{1}{p}}f^{*}(\theta x)\right\} =\frac{2^{\frac{1}{p}}B}{\theta^{\frac{1}{p}}}{\displaystyle \sup_{x>0}}\left\{ \left(\theta x\right)^{\frac{1}{p}}f^{*}(\theta x)\right\} =\frac{2^{\frac{1}{p}}B}{\theta^{\frac{1}{p}}}{\displaystyle \sup_{x>0}}\left\{ x^{\frac{1}{p}}f^{*}(x)\right\} .
\]
Letting $\theta\rightarrow1-$, it follows that $\left\Vert f\right\Vert _{L_{w\left(p,\infty\right)}^{\infty}}\leq2^{\frac{1}{p}}B\left\Vert f\right\Vert _{L\left(p,\infty\right)},$
and the last inequality in $(\ref{eq:estimate for GM for infinite q})$
follows.$\square$

The next lemma shows that the statement of Theorem 3.4 holds for $q<p$
if $f\in GM_{1}$. The proof is similar to the proof by B. Booton
{[}3{]} of the analogous statement for $GMS_{1}$.

$ $

\textbf{Lemma 4.9}. Let $f\in GM_{1}\left(B\right)$, $0<q<p$. Then
\[
\left\Vert f\right\Vert _{L(p,q)}\leq\left(\frac{2p}{q}\right)^{\frac{1}{q}}B^{2}\left\Vert f\right\Vert _{L_{w\left(p,q\right)}^{q}}.
\]

\textbf{Proof}.

Let $g(x)={\displaystyle \sum_{k=-\infty}^{\infty}}\left|f\left(2^{k}\right)\right|I_{\left[2^{k},2^{k+1}\right)}\left(x\right)$.
Then by $(\ref{eq:almost monotone-1})$, $\left|f\left(x\right)\right|\leq Bg\left(x\right)$.
Furthermore, $g^{*}$ is decreasing, right-continuous and for each
$x>0$, $g^{*}\left(x\right)=\left|f\left(2^{k}\right)\right|$ for
some $k\in\mathbb{Z}$. Therefore, $g^{*}$ can be written as 
\[
g^{*}\left(x\right){\displaystyle =\sum_{k\in J\cap\mathbb{Z}}}\alpha_{k}I_{\left[x_{k},x_{k+1}\right)}\left(x\right),
\]
where $J$ is an interval, $\alpha_{k}>a_{k+1}$, and $\left[x_{k},x_{k+1}\right)=\left\{ g^{*}=\alpha_{k}\right\} $,
$\forall k$. And so: 
\[
\left\Vert f\right\Vert _{L(p,q)}^{q}\leq B^{q}\left\Vert g\right\Vert _{L(p,q)}^{q}=B^{q}\intop_{0}^{\infty}x^{\frac{q}{p}}\left(g^{*}\left(x\right)\right)^{q}\frac{dx}{x}
\]
\[
=B^{q}\sum_{k=-\infty}^{\infty}\intop_{\left\{ g^{*}=\alpha_{k}\right\} }x^{\frac{q}{p}}\alpha_{k}^{q}\frac{dx}{x}=B^{q}\sum_{k=-\infty}^{\infty}\alpha_{k}^{q}\intop_{x_{k}}^{x_{k+1}}x^{\frac{q}{p}}\frac{dx}{x}
\]
\[
=B^{q}\frac{p}{q}\sum_{k=-\infty}^{\infty}\alpha_{k}^{q}\left(x_{k+1}^{\frac{q}{p}}-x_{k}^{\frac{q}{p}}\right)\leq B^{q}\frac{p}{q}\sum_{k=-\infty}^{\infty}\alpha_{k}^{q}\left(x_{k+1}-x_{k}\right)^{\frac{q}{p}}
\]
\[
=B^{q}\frac{p}{q}\sum_{k=-\infty}^{\infty}\alpha_{k}^{q}\left(\sum_{\left\{ m:\left|f\left(2^{m}\right)\right|=\alpha_{k}\right\} }2^{m}\right)^{\frac{q}{p}}\leq B^{q}\frac{p}{q}\sum_{k=-\infty}^{\infty}\sum_{\left\{ m:\left|f\left(2^{m}\right)\right|=\alpha_{k}\right\} }\left|f\left(2^{m}\right)\right|^{q}2^{\frac{mq}{p}}
\]
\[
=B^{q}\frac{p}{q}\sum_{m\in\mathbb{Z}}\left|f\left(2^{m}\right)\right|^{q}2^{\frac{mq}{p}}=2B^{q}\frac{p}{q}\sum_{m\in\mathbb{Z}}\intop_{2^{m-1}}^{2^{m}}\left|f\left(2^{m}\right)\right|^{\frac{q}{p}}2^{m\left(\frac{q}{p}-1\right)}dx
\]
\[
\leq2B^{2q}\frac{p}{q}\intop_{0}^{\infty}\left|f\left(x\right)\right|^{q}x^{\frac{q}{p}-1}dx.\square
\]

$ $

\textbf{Theorem 4.10}. Let $B>1$, $f\in GM_{1}\left(B\right)$ .
If $0<p<\infty$, $0<q<\infty$ then: 
\begin{equation}
\left\Vert f\right\Vert _{L(p,q)}\sim\left(\sum_{k=-\infty}^{\infty}2^{\frac{kq}{p}}\left(f^{*}\left(2^{k}\right)\right)^{q}\right)^{\frac{1}{q}}\sim\left(\sum_{k=-\infty}^{\infty}2^{\frac{kq}{p}}\left|f\left(2^{k}\right)\right|^{q}\right)^{\frac{1}{q}}\sim\left\Vert f\right\Vert _{L_{w(p,q)}^{q}}.\label{eq: equivalence of the two main norms on GM}
\end{equation}

If $0<p\leq\infty$ then
\begin{equation}
\left\Vert f\right\Vert _{L(p,\infty)}\sim\sup_{k\in\mathbb{Z}}{\scriptstyle \left\{ 2^{\frac{k}{p}}\left|f\left(2^{k}\right)\right|\right\} }\sim\sup_{k\in\mathbb{Z}}{\scriptstyle \left\{ 2^{\frac{k}{p}}\left|f^{*}\left(2^{k}\right)\right|\right\} \sim{\displaystyle \left\Vert f\right\Vert _{L_{w(p,\infty)}^{\infty}}}}.\label{eq: the missing link: case infinite q}
\end{equation}

\textbf{Proof}.

First, show that $\left\Vert f\right\Vert _{L(p,q)}\sim\left\Vert f\right\Vert _{L_{w(p,q)}^{q}}$
for all $0<q\leq\infty$.

The inequality $\left\Vert f\right\Vert _{L_{w\left(p,q\right)}^{q}}\leq B^{2}C(p,q)\left\Vert f\right\Vert _{L\left(p,q\right)}$
follows from Lemma 4.8, more specifically, for $0<p<\infty$, $0<q<\infty$
it is a consequence of $(\ref{eq: another estimate between Lqw and})$,
and for $0<p\leq\infty$, $q=\infty$ it is a consequence of $(\ref{eq:estimate for GM for infinite q})$.

The inequality $\left\Vert f\right\Vert _{L(p,q)}\leq B^{2}C(p,q)\left\Vert f\right\Vert _{L_{w\left(p,q\right)}^{q}}$
follows from Theorem 3.4 for all $0<p\leq q\leq\infty$ and it is
a consequence of Lemma 4.9 for all $0<q<p$.

$f^{*}\in GMS_{1}\left(1\right)$ as a decreasing function, and so
by Lemma 4.8, if $0<p<\infty$ then $\left({\displaystyle \sum_{k=-\infty}^{\infty}}2^{\frac{kq}{p}}\left|f^{*}\left(2^{k}\right)\right|^{q}\right)^{\frac{1}{q}}\sim\left\Vert f\right\Vert _{L\left(p,q\right)}$,
and ${\displaystyle \sup_{k\in\mathbb{Z}}}\left\{ 2^{\frac{k}{p}}\left|f^{*}\left(2^{k}\right)\right|\right\} \sim\left\Vert f\right\Vert _{L\left(p,\infty\right)}$,
that is, the remaining equivalences in $(\ref{eq: equivalence of the two main norms on GM})$
and $(\ref{eq: the missing link: case infinite q})$ follow$.\square$

$ $

\textbf{5. Properties of $GM_{\alpha,\varphi}$ and $GMS_{\alpha,\varphi}$.}

$ $

\textbf{Definition 5.1}. For each $f\in L_{loc}^{1}\left(\mathbb{R}^{+},\mu\right)$,
where $\mu$ is a Borel measure, define for each $x$ such that $\mu\left(\left(0,x\right]\right)>0$,
the \emph{average of $f$ on $\left(0,x\right]$, }as 
\begin{equation}
\sigma_{x}(f,\mu):={\displaystyle \frac{1}{\mu\left((0,x]\right)}\intop_{(0,x]}}fd\mu.\label{eq: average of f}
\end{equation}

Formula $(\ref{eq: average of f})$ defines a function of $x$, which
we shall write as $\sigma(f,\mu)$. If $f$ or $\mu$ are clear from
the context, we may omit them from the notation.

$ $

\textbf{Lemma 5.2}. Let $(\Omega,\Sigma,\mu)$ be a measure space,
$-\infty<c<d<\infty$, and assume that $g:\Omega\times\mathbb{R}\rightarrow\mathbb{C}$
is such that for all $c\leq x\leq d$, $g(\cdot,x)$ is a $\mu$-integrable
function. Then: 
\begin{equation}
V_{\intop_{\Omega}g(\omega,\cdot)d\mu}([c,d])\leq\intop_{\Omega}V_{g(\omega,\cdot)}([c,d])d\mu.\label{eq:var of int is < int of var}
\end{equation}

\textbf{Lemma 5.3}. Let $(\Omega,\Sigma,\mu)$ be a measure space,
and assume that for a.e. $\omega\in\Omega$, $g(\omega,\cdot)\in GM_{\alpha,\varphi}(B)$.
Then 
\begin{equation}
f(\cdot):=\intop_{\Omega}g(\omega,\cdot)d\mu\in GM_{\alpha,\varphi}\left(\frac{B}{\cos\varphi}\right).\label{eq:GM is an integral cone}
\end{equation}

\textbf{Corollary 5.3}. $GM_{\alpha,\varphi}$ is a cone in $GM$.
Moreover, if $f_{j}\in GM_{\alpha,\varphi}(B)$, $1\leq j\leq N$
then: 
\begin{equation}
{\displaystyle \sum_{j=1}^{N}}f_{j}\in GM_{\alpha,\varphi}\left(\frac{B}{\cos\varphi}\right).\label{eq:GM(a,phi)  is a cone}
\end{equation}

Similarly, $GMS_{\alpha,\varphi}$ is a cone in $GMS$ and for $b_{n}={\displaystyle \sum_{j=1}^{N}}a_{n,j}$,
where for each $j$ the sequence $\left\{ a_{n,j}\right\} _{n=1}^{\infty}$
is in $GMS(B)$, we have that 
\begin{equation}
\left\{ b_{n}\right\} _{n=1}^{\infty}\in GMS_{\alpha,\varphi}\left(\frac{B}{\cos\varphi}\right).\square\label{eq:GMS(a,phi) is a cone}
\end{equation}

\textbf{Lemma 5.4}. Assume that for a.e. $x>0$, $f(x)\in S_{\alpha,\varphi}$,
and that $\mu$ is a Borel measure on $\mathbb{R}^{+}$ for which
there is $A\geq0$ such that for all $a>A$, $\mu\left(\left(0,a\right]\right)>0$.
Then, $\forall x>A$: 
\begin{equation}
\sigma_{x}(\left|f\right|)\leq\frac{1}{\cos\varphi}\left|\sigma_{x}(f)\right|.\label{eq:averages in the cones: f and |f|}
\end{equation}

If, moreover, $f\in GM_{1}\left(B\right)$ and $\exists K>0$ such
that $\forall a>0$, 
\begin{equation}
\mu((0,a])\leq K\mu((a,2a])\label{eq: admissibilty of the measure for averages-1}
\end{equation}
then, $\forall x>A$: 
\begin{equation}
\left|f(x)\right|\leq\frac{\left(K+1\right)B}{\cos\varphi}\left|\sigma_{x}(f)\right|.\label{eq:averages in the cones}
\end{equation}

\textbf{Proof}.

$(\ref{eq:averages in the cones: f and |f|})$ follows from Lemma
1.19:
\[
\sigma_{x}(\left|f\right|)={\displaystyle \frac{1}{\mu\left((0,x]\right)}\intop_{(0,x]}}\left|f\right|d\mu\leq\frac{1}{\mu\left((0,x]\right)\cos\varphi}\left|\intop_{(0,x]}fd\mu\right|=\frac{1}{\cos\varphi}\left|\sigma_{x}(f)\right|.
\]
\[
\]
Since $f\in GM_{1}\left(B\right)$, it follows that $\left|f(x)\right|\leq B\left|f(t)\right|$,
$\forall t\in\left[{\displaystyle \frac{x}{2}},x\right]$. Furthermore,
by $(\ref{eq: admissibilty of the measure for averages-1})$ for all
$x>0$: 
\[
\mu\left((0,x]\right)=\mu\left(\left(0,{\textstyle \frac{x}{2}}\right]\right)+\mu\left(\left({\textstyle \frac{x}{2}},x\right]\right)\leq\left(K+1\right)\mu\left(\left({\textstyle \frac{x}{2}},x\right]\right).
\]
Therefore, for $x>A$: 
\[
\left|f(x)\right|\leq\frac{B}{\mu\left(\left(\frac{x}{2},x\right]\right)}\intop_{\left({\textstyle \frac{x}{2}},x\right]}\left|f(t)\right|d\mu\leq\frac{\left(K+1\right)B}{\mu\left((0,x]\right)}\intop_{(0,x]}\left|f(t)\right|d\mu
\]
\[
=\left(K+1\right)B\sigma_{x}(\left|f\right|)\leq\frac{\left(K+1\right)B}{\cos\varphi}\left|\sigma_{x}(f)\right|.\square
\]

Analogously,

\textbf{Lemma 5.5}. Let $a_{k}\in S_{\alpha,\varphi}$, $\forall k\geq1$.
Then, $\forall n\geq1$: 
\[
\sigma_{n}\left(\left\{ \left|a_{k}\right|\right\} \right)\leq\frac{1}{\cos\varphi}\left|\sigma_{n}\left(\left\{ a_{k}\right\} \right)\right|.
\]

If, moreover, $\left\{ a_{k}\right\} \in GMS_{1}\left(B\right)$,
then, $\forall n\geq1$:

\[
\left|a_{n}\right|\leq\frac{2B}{\cos\varphi}\left|\sigma_{n}\left(\left\{ a_{k}\right\} \right)\right|.
\]
$ $

\textbf{Lemma 5.6.} Let $f\in GM_{\alpha,\varphi}\left(B\right)$
and assume that $\mu$ is a Borel measure on $\mathbb{R}^{+}$, for
which there are $A\geq0$, $K>0$ such that for all $a>A$, 
\begin{equation}
\begin{cases}
\begin{array}{c}
\mu\left(\left(0,a\right]\right)>0\\
\mu((a,2a])\leq K\mu((0,a])\leq K^{2}\mu\left(\left(0,a\right]\right)
\end{array} & .\end{cases}\label{eq: admissibilty of the measure for averages}
\end{equation}

Then $\forall x>A$, $\sigma_{x}(f,\mu)\in S_{\alpha,\varphi}$. Moreover,
there is $B'=B'\left(B,K,\varphi\right)>0$ such that for all $x>A$,
$\sigma(f,\mu)$ satisfies 
\begin{equation}
V_{\sigma(f,\mu)}\left(\left[x,2x\right]\right)\leq B'\left|\sigma_{x}(f,\mu)\right|.\label{eq:GM with delay}
\end{equation}

In particular, if $A=0$ then $\sigma(f,\mu)\in GM_{\alpha,\varphi}\left(B'\right)$.

\textbf{Proof}.

Clearly, $\sigma_{x}(f,\mu)\in S_{\alpha,\varphi}$.

Let $a>A$. Define, for $x\geq a$: 
\begin{equation}
g(t,x)=f(t)\frac{I_{(0,x]}(t)}{\mu((0,x])}.\label{eq:g(omega,x)}
\end{equation}

Since $f\in L_{loc}^{1}\left(\mathbb{R}^{+},\mu\right)$, it follows
that $g(\cdot,x)\in L^{1}\left(\mathbb{R}^{+},\mu\right)$. Furthermore,
$\sigma_{x}(f)={\displaystyle \intop_{\mathbb{R}^{+}}g(t,x)d\mu}$.
By Lemma 5.2:
\[
V_{\sigma(f)}\left([a,2a]\right)=V_{\underset{\mathbb{R}^{+}}{\int}g(t,\cdot)d\mu}\left([a,2a]\right)\leq{\displaystyle \intop_{\mathbb{R}^{+}}}V_{g(t,\cdot)}\left([a,2a]\right)d\mu.
\]

Clearly, $g\left(t,x\right)=f(t){\displaystyle \frac{I_{[t,\infty)}(x)}{\mu((0,x])}}$
is a decreasing function of $x$ on $\left[t,\infty\right)$. Consider
the following three cases: $0\leq t\leq a$, $a<t\leq2a$, $t>2a$.

If $0\leq t\leq a$ then 
\[
V_{g\left(t,x\right)}\left(\left[a,2a\right]\right)=V_{\frac{{\scriptstyle f\left(t\right)}}{{\scriptstyle \mu\left(\left(0,x\right]\right)}}}\left(\left[a,2a\right]\right)=\frac{\left|f\left(t\right)\right|}{\mu\left(\left(0,a\right]\right)}-\frac{\left|f\left(t\right)\right|}{\mu\left(\left(0,2a\right]\right)}.
\]

If $a<t\leq2a$, observe that for $a\leq x<t$, $g\left(t,x\right)=0$,
and $g(t,t)={\displaystyle \frac{\left|f\left(t\right)\right|}{\mu\left(\left(0,t\right]\right)}}$,
therefore, $V_{g\left(t,x\right)}\left(\left[a,t\right]\right)={\displaystyle \frac{\left|f\left(t\right)\right|}{\mu\left(\left(0,t\right]\right)}}$,
and so: 
\[
V_{g\left(t,x\right)}\left(\left[a,2a\right]\right)=V_{g\left(t,x\right)}\left(\left[a,t\right]\right)+V_{g\left(t,x\right)}\left(\left[t,2a\right]\right)
\]
\[
=\frac{\left|f\left(t\right)\right|}{\mu\left(\left(0,t\right]\right)}+V_{\frac{{\scriptstyle f\left(t\right)}}{{\scriptstyle \mu\left(\left(0,x\right]\right)}}}\left(\left[t,2a\right]\right)
\]
\[
=\frac{\left|f\left(t\right)\right|}{\mu\left(\left(0,a\right]\right)}+\left(\frac{\left|f\left(t\right)\right|}{\mu\left(\left(0,t\right]\right)}-\frac{\left|f\left(t\right)\right|}{\mu\left(\left(0,2a\right]\right)}\right)=\frac{2\left|f\left(t\right)\right|}{\mu\left(\left(0,t\right]\right)}-\frac{\left|f\left(t\right)\right|}{\mu\left(\left(0,2a\right]\right)}.
\]

Finally, if $t>2a$ then $g\left(t,x\right)=0$, $\forall x\in\left[a,2a\right]$,
and so, $V_{g\left(t,x\right)}\left(\left[a,2a\right]\right)=0$.
Therefore: 
\[
{\displaystyle \intop_{\mathbb{R}^{+}}}V_{g(t,\cdot)}\left([a,2a]\right)d\mu=
\]
\[
=\intop_{(0,a]}\left(\frac{1}{\mu((0,a])}-\frac{1}{\mu((0,2a])}\right)\left|f(t)\right|d\mu
\]
\[
+\intop_{(a,2a]}\left(\frac{2}{\mu((0,t])}-\frac{1}{\mu((0,2a])}\right)\left|f(t)\right|d\mu
\]
\[
\leq\left(\frac{2}{\mu((0,a])}-\frac{1}{\mu((0,2a])}\right)\intop_{(0,2a]}\left|f(t)\right|d\mu
\]
\[
=\frac{2\mu((a,2a])+\mu\left(\left(0,a\right]\right)}{\mu((0,a])}\sigma_{2a}(\left|f\right|)\leq\left(2K+1\right)\sigma_{2a}(\left|f\right|).
\]

Observe that by Lemma 1.22 for $f\in GM\left(B\right)$ and $a<t\leq2a$,
$\left|f\left(t\right)\right|\leq2B\left|f\left(a\right)\right|$,
and so: 
\[
\sigma_{2a}(\left|f\right|)=\frac{1}{\mu\left(\left(0,2a\right]\right)}\left(\intop_{\left(0,a\right]}\left|f\right|d\mu+\intop_{\left(a,2a\right]}\left|f\right|d\mu\right)
\]
\[
\leq\frac{1}{\mu\left(\left(0,a\right]\right)}\intop_{\left(0,a\right]}\left|f\right|d\mu+\frac{1}{\mu\left(\left(a,2a\right]\right)}\intop_{\left(a,2a\right]}2B\left|f\left(a\right)\right|d\mu=\sigma_{a}(\left|f\right|)+2B\left|f\left(a\right)\right|.
\]

Since $\mu$ satisfies $(\ref{eq: admissibilty of the measure for averages-1})$
and $f\in GM_{\alpha,\varphi}$, it follows from Lemma 5.4 that:

\[
\sigma_{a}(\left|f\right|)\leq\frac{1}{\cos\varphi}\left|\sigma_{a}(f)\right|\mbox{ and }\left|f(a)\right|\leq\frac{\left(K+1\right)B}{\cos\varphi}\left|\sigma_{a}(f)\right|,
\]
and so, $\sigma_{2a}(\left|f\right|)\leq{\displaystyle \frac{1+2\left(K+1\right)B^{2}}{\cos\varphi}}\left|\sigma_{a}(f)\right|$.
Therefore,

\[
V_{\sigma(f,\mu)}\left([a,2a]\right)\leq\left(2K+1\right)\sigma_{2a}(\left|f\right|)\leq\frac{\left(2K+1\right)\left(1+2\left(K+1\right)B^{2}\right)}{\cos\varphi}\left|\sigma_{a}(f,\mu)\right|.\square
\]

$ $

\textbf{6.} \textbf{Fourier series with general monotone coefficients.}

$ $

\textbf{Theorem 6.1}. Assume that $a_{k}\in\mathbb{C}$ , $1\leq m\leq N$.
Then, $\forall x\in(0,\pi]$: 
\begin{equation}
\left|\sum_{k=m}^{N}a_{k}e^{ikx}\right|\leq\frac{4\pi}{x}\left(\frac{\left|a_{m}\right|}{2}+\sum_{k=m}^{N-1}\left|a_{k+1}-a_{k}\right|\right).\label{eq:boundary by the variation of the an. F.s.}
\end{equation}

\textbf{Proof.}

Since for $x\in\left(0,\pi\right]$, 

\begin{equation}
\left|\sum_{k=0}^{N}e^{ikx}\right|\leq\frac{1}{\sin\frac{x}{2}}\leq\frac{\pi}{x},\label{eq:est on th Conj. Poisson Kernel}
\end{equation}

it follows, using summation by parts,
\[
\sum_{k=m}^{N}a_{k}e^{ikx}=a_{N}\sum_{k=m}^{N}e^{ikx}-\sum_{k=m}^{N-1}\left(a_{k+1}-a_{k}\right)\sum_{j=m}^{k}e^{ijx},
\]

and so 
\[
\left|\sum_{k=m}^{N}a_{k}e^{ikx}\right|\leq\frac{2\pi}{x}\left(\left|a_{N}\right|+\sum_{k=m}^{N-1}\left|a_{k+1}-a_{k}\right|\right).
\]

Since $a_{N}=a_{m}+{\displaystyle \sum_{k=m}^{N-1}}\left(a_{k+1}-a_{k}\right)$,
it follows that 
\begin{equation}
\left|\sum_{k=m}^{N}a_{k}e^{ikx}\right|\leq\frac{2\pi}{x}\left(\left|a_{m}\right|+2\cdot\sum_{k=m}^{N-1}\left|a_{k+1}-a_{k}\right|\right).\square\label{eq:variation of sine F.S.}
\end{equation}

\textbf{Corollary 6.1}. Assume $\left\{ a_{k}\right\} \in GMS_{2}\left(B\right)$.
Then $\forall x\in(0,\pi]$, $\forall1\leq m\leq N$, 
\begin{equation}
\left|\sum_{k=m}^{N}a_{k}e^{ikx}\right|\leq\frac{6\pi B}{x}\left(\left|a_{m}\right|+\sum_{k=m+1}^{N}\frac{\left|a_{k}\right|}{k}\right).\label{eq:first estimate for GM series}
\end{equation}

$ $

\textbf{Theorem 6.2}. Assume $\left\{ c_{k}\right\} \in GMS_{2}\left(B\right)$.
If $f(x)={\displaystyle \sum_{k=1}^{\infty}}c_{k}e^{ikx},$ then:

\begin{equation}
\left\Vert f\right\Vert _{L^{1}(0,\pi)}\leq2\pi\left|c_{1}\right|+27\pi B\left\Vert \left\{ c_{n}\right\} \right\Vert _{l_{\frac{\ln k}{k}}^{1}}.\label{eq:L1 of GM Fourier series}
\end{equation}

\textbf{Proof.}

For each $N\geq2$, denote $f_{N}(x)={\displaystyle \sum_{k=1}^{N}}c_{k}e^{ikx}$.
Then:

\[
\left\Vert f_{N}\right\Vert _{L^{1}(0,\pi)}=\sum_{n=1}^{\infty}\intop_{\frac{\pi}{n+1}}^{\frac{\pi}{n}}\left|{\displaystyle \sum_{k=1}^{N}}c_{k}e^{ikx}\right|dx
\]
\[
\leq\sum_{n=1}^{N-1}\intop_{\frac{\pi}{n+1}}^{\frac{\pi}{n}}\left|{\displaystyle \sum_{k=1}^{N}}c_{k}e^{ikx}\right|dx+\sum_{n=N}^{\infty}\intop_{\frac{\pi}{n+1}}^{\frac{\pi}{n}}\left({\displaystyle \sum_{k=1}^{N}}\left|c_{k}\right|\right)dx
\]
\[
\leq\sum_{n=1}^{N-1}\intop_{\frac{\pi}{n+1}}^{\frac{\pi}{n}}\left(\sum_{k=1}^{n}\left|c_{k}\right|+\left|\sum_{k=n+1}^{N}c_{k}e^{ikx}\right|\right)dx+\frac{\pi}{N}\sum_{k=1}^{N}\left|c_{k}\right|.
\]

Applying $\ensuremath{(\ref{eq:first estimate for GM series})}$ for
$n<N-1$: 
\[
\left|\sum_{k=n+1}^{N}c_{k}e^{ikx}\right|\leq\frac{6\pi B}{x}\left(\left|c_{n+1}\right|+\sum_{k=n+2}^{N}\frac{\left|c_{k}\right|}{k}\right).
\]
If $n=N-1$ then $\left|{\displaystyle \sum_{k=n+1}^{N}}c_{k}e^{ikx}\right|=\left|c_{N}\right|$,
and so, for $n\leq N-1$,
\[
\left|\sum_{k=n+1}^{N}c_{k}e^{ikx}\right|\leq\frac{6\pi B}{x}\left(\left|c_{n+1}\right|+\sum_{k=n+1}^{N}\frac{\left|c_{k}\right|}{k}\right).
\]
Therefore:
\[
\left\Vert f_{N}\right\Vert _{L^{1}(0,\pi)}\leq\sum_{n=1}^{N-1}\intop_{\frac{\pi}{n+1}}^{\frac{\pi}{n}}\left[\sum_{k=1}^{n}\left|c_{k}\right|+\frac{6\pi B}{x}\left(\left|c_{n+1}\right|+\sum_{k=n+1}^{N}\frac{\left|c_{k}\right|}{k}\right)\right]dx+\frac{\pi}{N}\sum_{k=1}^{N}\left|c_{k}\right|
\]
\[
=I_{1}+I_{2}+I_{3}+\frac{\pi}{N}\sum_{k=1}^{N}\left|c_{k}\right|,
\]
where, using summation by parts: 
\[
I_{1}=\sum_{n=1}^{N-1}\intop_{\frac{\pi}{n+1}}^{\frac{\pi}{n}}\left(\sum_{k=1}^{n}\left|c_{k}\right|\right)dx=\sum_{n=1}^{N-1}\left(\frac{\pi}{n}-\frac{\pi}{n+1}\right)\sum_{k=1}^{n}\left|c_{k}\right|
\]

\[
=-\frac{\pi}{N}\sum_{n=1}^{N}\left|c_{n}\right|+\pi\sum_{n=1}^{N}\frac{\left|c_{n}\right|}{n}\leq{\displaystyle \pi\sum_{n=1}^{N}\frac{\left|c_{n}\right|}{n}}.
\]
It follows that 
\[
I_{1}\leq{\displaystyle \pi\left(\left|c_{1}\right|+\sum_{n=2}^{N}\frac{\left|c_{n}\right|}{n}\right){\displaystyle \leq\pi\left|c_{1}\right|+\frac{\pi}{\ln2}\sum_{n=2}^{N}\left|c_{n}\right|\frac{\ln n}{n}}}.
\]
\[
I_{2}=\sum_{n=1}^{N-1}\intop_{\frac{\pi}{n+1}}^{\frac{\pi}{n}}\frac{6\pi B}{x}\left|c_{n+1}\right|dx=
\]
\[
=6\pi B\sum_{n=1}^{N-1}\left|c_{n+1}\right|\ln\frac{n+1}{n}\leq6\pi B\sum_{n=1}^{N-1}\frac{\left|c_{n+1}\right|}{n}=6\pi B\sum_{n=2}^{N}\frac{\left|c_{n}\right|}{n-1}
\]
\[
\leq12\pi B\sum_{n=2}^{N}\frac{\left|c_{n}\right|}{n}\leq\frac{12\pi B}{\ln2}\sum_{n=2}^{N}\left|c_{n}\right|\frac{\ln n}{n};
\]
\[
I_{3}=\sum_{n=1}^{N-1}\intop_{\frac{\pi}{n+1}}^{\frac{\pi}{n}}\frac{6\pi B}{x}\left(\sum_{k=n+1}^{N}\frac{\left|c_{k}\right|}{k}\right)dx=6\pi B\sum_{n=1}^{N-1}\left(\sum_{k=n+1}^{N}\frac{\left|c_{k}\right|}{k}\right)\intop_{\frac{\pi}{n+1}}^{\frac{\pi}{n}}\frac{dx}{x}
\]
\[
=6\pi B\sum_{n=1}^{N-1}\sum_{k=n}^{N-1}\frac{\left|c_{k+1}\right|}{k+1}\ln\frac{n+1}{n}=6\pi B\sum_{k=1}^{N-1}\frac{\left|c_{k+1}\right|}{k+1}\sum_{n=1}^{k}\ln\frac{n+1}{n}
\]
\[
=6\pi B\sum_{k=1}^{N-1}\frac{\left|c_{k+1}\right|}{k+1}\ln\left(k+1\right)=6\pi B\sum_{k=2}^{N}\left|c_{k}\right|\frac{\ln k}{k}.
\]

Combining the estimates for $I_{j}$ and taking into account that
\[
{\displaystyle \frac{\pi}{N}\sum_{k=1}^{N}}\left|c_{k}\right|\leq\pi\left|c_{1}\right|+\frac{\pi}{N\ln2}{\displaystyle \sum_{k=2}^{N}}\left|c_{k}\right|\ln k\leq\pi\left|c_{1}\right|+\frac{\pi B}{\ln2}{\displaystyle \sum_{k=2}^{N}}\left|c_{k}\right|{\displaystyle \frac{\ln k}{k}},
\]
we have: 
\[
\left\Vert f_{N}\right\Vert _{L^{1}(0,\pi)}\leq2\pi\left|c_{1}\right|+27\pi B\sum_{k=2}^{N}\left|c_{k}\right|\frac{\ln k}{k}.\square
\]

$ $

\textbf{Theorem 6.3}. Assume $\left\{ c_{k}\right\} \in GMS_{2}\left(B\right)$.
If $f(x)={\displaystyle \sum_{k=1}^{\infty}}c_{k}e^{ikx}$, then

\begin{equation}
\left\Vert f\right\Vert _{L(1,\infty)(0,\pi)}\leq6\pi B\left\Vert \left\{ c_{n}\right\} \right\Vert _{l_{1/k}^{1}}.\label{eq:weak-L1 of GM2 Fourier series}
\end{equation}

\textbf{Proof}.

Let $N\geq1$. Denote $f_{N}(x)={\displaystyle \sum_{k=1}^{N}}c_{k}e^{ikx}$.
Then, by $(\ref{eq:first estimate for GM series})$, 
\[
\left|f_{N}(x)\right|=\left|\sum_{k=1}^{N}c_{k}e^{ikx}\right|\leq\frac{6\pi B}{x}\left(\left|c_{1}\right|+\sum_{k=2}^{N}\frac{\left|c_{k}\right|}{k}\right)=\frac{6\pi B}{x}\sum_{k=1}^{N}\frac{\left|c_{k}\right|}{k}.
\]

Therefore, for $\forall\alpha>0$: 
\[
\lambda\left(\left\{ x:\left|f_{N}(x)\right|>\alpha\right\} \right)\leq\lambda\left(\left\{ x:\frac{6\pi B}{x}\sum_{k=1}^{N}\frac{\left|c_{k}\right|}{k}>\alpha\right\} \right)
\]
\[
\leq\lambda\left(\left\{ x:\frac{6\pi B}{x}\left\Vert \left\{ c_{k}\right\} \right\Vert _{l_{\frac{1}{k}}^{1}}>\alpha\right\} \right)=\frac{6\pi B}{\alpha}\left\Vert \left\{ c_{k}\right\} \right\Vert _{l_{\frac{1}{k}}^{1}},
\]

proving $(\ref{eq:weak-L1 of GM2 Fourier series})$ $\square$

The result of Theorem 6.3 can be interpreted as follows. Let $T$
be the operator mapping $GMS_{2}$ to functions over $\mathbb{R}$,
defined as 
\begin{equation}
T\left(\left\{ c_{n}\right\} \right)(x):=\sum_{k=1}^{\infty}c_{k}e^{ikx}.\label{eq:Fourier series operator on GM - definition}
\end{equation}

Then 
\begin{equation}
T:l_{\frac{1}{k}}^{1}\cap GMS_{2}\rightarrow L(1,\infty).\label{eq:Fourier series operator on GM}
\end{equation}

$ $

\textbf{Theorem 6.4}. Let $1<p<\infty$, $0<q\leq\infty$. If $c_{k}$
are trigonometric Fourier coefficients of $f$ and $\left\{ c_{k}\right\} \in GMS_{\alpha,\varphi}(B)$
then $\exists C=C(p,q,\varphi)$ such that 
\begin{equation}
\left\Vert \left\{ c_{n}\right\} \right\Vert _{l(p',q)}\leq C\left\Vert f\right\Vert _{L(p,q)(0,2\pi)}.\label{eq:Main result-1}
\end{equation}

\textbf{Proof.}

By Lemma 5.5, $\left|c_{n}\right|\leq{\displaystyle \frac{2B}{\cos\varphi}}\left|\sigma_{n}\right|$,
and so $\left\Vert \left\{ c_{n}\right\} \right\Vert _{l(p',q)}\leq{\displaystyle \frac{2B}{\cos\varphi}}\left\Vert \sigma_{n}\right\Vert _{l(p',q)}$.
By $(\ref{eq: Sagher-12-1})$, 
\[
\]
\[
\left\Vert \left\{ c_{n}\right\} \right\Vert _{l(p',q)}\leq\frac{C'\left(p,q\right)}{\cos\varphi}\left\Vert f\right\Vert _{L(p,q)(0,2\pi)}.\square
\]

$ $

\textbf{7. Interpolation in $GMS_{\alpha,\varphi}$.}

$ $

The following Lemma is a special case of a result of Yu. Brudnui and
N. Krugljak ({[}4{]}, Corollary 3.1.26).

$ $

\textbf{Lemma 7.1}. 
\begin{equation}
K\left(t,\left\{ c_{n}\right\} ,l_{\frac{1}{k}}^{1},l^{1}\right)\sim t\sum_{n\leq\frac{1}{t}}\left|c_{n}\right|+\sum_{n>\frac{1}{t}}\frac{\left|c_{n}\right|}{n}\label{eq:K between lw and l}
\end{equation}

$ $

\textbf{Lemma 7.2}. Let $\left\{ a_{n}\right\} ,\left\{ c_{n}\right\} \in GMS\left(B\right)$.
Assume that $N\geq1$ and for some $\gamma>0$, $\left|c_{N}\right|\leq\gamma\left|a_{N}\right|$.
Define 
\[
b_{n}=\begin{cases}
\begin{array}{cc}
a_{n} & \mbox{if }n\leq N\\
c_{n} & \mbox{if }n>N
\end{array} & .\end{cases}
\]
Then $\left\{ b_{n}\right\} \in GMS\left(3B+6B^{2}\gamma\right)$.

\textbf{Proof}.

If $n=N=1$ then 
\[
{\displaystyle \sum_{k=n}^{2n-1}}\left|b_{k+1}-b_{k}\right|=\left|a_{1}-c_{2}\right|=\left|a_{N}-c_{N+1}\right|.
\]

If $n=N>1$ then 
\[
{\displaystyle \sum_{k=n}^{2n-1}}\left|b_{k+1}-b_{k}\right|=\left|a_{N}-c_{N+1}\right|+{\displaystyle \sum_{k=N+1}^{2N-1}}\left|c_{k+1}-c_{k}\right|.
\]

If $n={\displaystyle \frac{N+1}{2}}>1$ then 
\[
{\displaystyle \sum_{k=n}^{2n-1}}\left|b_{k+1}-b_{k}\right|={\displaystyle \sum_{k=n}^{N-1}}\left|a_{k+1}-a_{k}\right|+\left|a_{N}-c_{N+1}\right|.
\]

If ${\displaystyle \frac{N+1}{2}}<n<N$ then 
\[
{\displaystyle \sum_{k=n}^{2n-1}}\left|b_{k+1}-b_{k}\right|={\displaystyle \sum_{k=n}^{N-1}}\left|a_{k+1}-a_{k}\right|+\left|a_{N}-c_{N+1}\right|+{\displaystyle \sum_{k=N+1}^{2n-1}}\left|c_{k+1}-c_{k}\right|.
\]

In all four cases above: 
\[
{\displaystyle \sum_{k=n}^{2n-1}}\left|b_{k+1}-b_{k}\right|\leq{\displaystyle \sum_{k=n}^{2n-1}}\left|a_{k+1}-a_{k}\right|+\left|a_{N}-c_{N+1}\right|+{\displaystyle \sum_{k=N}^{2N-1}}\left|c_{k+1}-c_{k}\right|,
\]

and ${\displaystyle \frac{N}{2}}<n\leq N$. By $(\ref{eq:almost monotone: B})$,
$\left|a_{N}\right|\leq2B\left|a_{n}\right|$ and $\left|c_{N+1}\right|\leq2B\left|c_{N}\right|$,
and so: 
\[
{\displaystyle \sum_{k=n}^{2n-1}}\left|b_{k+1}-b_{k}\right|\leq{\displaystyle \sum_{k=n}^{2n-1}}\left|a_{k+1}-a_{k}\right|+\left|a_{N}\right|+\left|c_{N+1}\right|+{\displaystyle \sum_{k=N}^{2N-1}}\left|c_{k+1}-c_{k}\right|
\]
\[
\leq B\left|a_{n}\right|+2B\left|a_{n}\right|+(2B+B)\left|c_{N}\right|\leq3B\left|a_{n}\right|+3B\gamma\left|a_{N}\right|
\]
\[
\leq3B\left|a_{n}\right|+6B^{2}\gamma\left|a_{n}\right|=\left(3B+6B^{2}\gamma\right)\left|b_{n}\right|.
\]

Two cases remain. If $n<{\displaystyle \frac{N+1}{2}}$ then $2n-1\leq N-1$,
and so: 
\[
{\displaystyle \sum_{k=n}^{2n-1}}\left|b_{k+1}-b_{k}\right|={\displaystyle \sum_{k=n}^{2n-1}}\left|a_{k+1}-a_{k}\right|\leq B\left|a_{n}\right|=B\left|b_{n}\right|.
\]

If $n>N$ then 
\[
{\displaystyle \sum_{k=n}^{2n-1}}\left|b_{k+1}-b_{k}\right|={\displaystyle \sum_{k=n}^{2n-1}}\left|c_{k+1}-c_{k}\right|\leq B\left|c_{n}\right|=B\left|b_{n}\right|.\square
\]

$ $

Lemmas 7.1 and 7.2 imply:

$ $

\textbf{Theorem 7.3}. 
\begin{equation}
K\left(t,\left\{ c_{n}\right\} ,l_{\frac{1}{k}}^{1}\cap GMS_{\alpha,\varphi},l^{1}\right)\leq\frac{9}{2}K\left(t,\left\{ c_{n}\right\} ,l_{\frac{1}{k}}^{1},l^{1}\right).\label{eq:K on GMS}
\end{equation}

\textbf{Proof}.

Let $\left\{ c_{n}\right\} \in GMS_{\alpha,\varphi}(B)\cap\left(l_{\frac{1}{k}}^{1}+l^{1}\right)=GMS_{\alpha,\varphi}(B)\cap l_{\frac{1}{k}}^{1}$.

1) If $t>1$ then, by $(\ref{eq:K between lw and l})$, $K\left(t,\left\{ c_{n}\right\} ,l_{\frac{1}{k}}^{1},l^{1}\right)={\displaystyle \sum_{n=1}^{\infty}\frac{\left|c_{n}\right|}{n}}=\left\Vert \left\{ c_{n}\right\} \right\Vert _{l_{\frac{1}{k}}^{1}}$,
and $(\ref{eq:K on GMS})$ follows.

2) Assume that $t\leq1$ .

Let $N=1+{\displaystyle \left\lfloor \frac{1}{t}\right\rfloor }$,
$\sigma_{N}={\displaystyle \frac{1}{N}}{\displaystyle \sum_{k=1}^{N}}\left|c_{k}\right|$,
and define $a_{n}={\displaystyle \frac{n}{N}}\sigma_{N}e^{i\alpha}$
for all $n\geq1$. Then ${\displaystyle \sum_{k=n}^{2n-1}}\left|a_{k}-a_{k+1}\right|=\left|a_{2n}-a_{n}\right|={\displaystyle \frac{n}{N}}\sigma_{N}=\left|a_{n}\right|$,
and so $\left\{ a_{n}\right\} \in GMS(1)$. By Corollary 5.3, $\left\{ a_{n}\right\} \in GMS_{\alpha,\varphi}(1)$.
Observe that for $2\leq N\leq4$, $\#\left\{ k\in\mathbb{Z}:\frac{N}{4}<k\leq\frac{N}{2}\right\} =1>{\displaystyle \frac{N}{5}}$.
Assume that $N>4$ and let $\theta=\frac{N}{4}-\left\lfloor \frac{N}{4}\right\rfloor $.
If $\theta<\frac{1}{2}$ then $\theta\leq\frac{1}{4}$, and so: 
\[
\#\left\{ k\in\mathbb{Z}:\frac{N}{4}<k\leq\frac{N}{2}\right\} =\left\lfloor \frac{N}{4}\right\rfloor \geq\frac{N}{4}-\frac{1}{4}\geq\frac{N}{4}-\frac{N/5}{4}=\frac{N}{5}.
\]

If $\theta\geq\frac{1}{2}$ then 
\[
\#\left\{ k\in\mathbb{Z}:\frac{N}{4}<k\leq\frac{N}{2}\right\} \geq\left\lfloor \frac{N}{4}\right\rfloor +1>\frac{N}{4}>{\displaystyle \frac{N}{5}.}
\]

Applying $(\ref{eq:almost monotone-1})$ to $\left\{ c_{n}\right\} $,
for all $n\in\left[{\displaystyle \frac{N}{2}},N\right]$ we have:
\[
\left|a_{n}\right|=\frac{n}{N}\sigma_{N}\geq\frac{1}{2}\sigma_{N}\geq\frac{1}{2}\cdot\frac{1}{N}\sum_{k\in\left(\frac{N}{4},\frac{N}{2}\right]}\left|c_{k}\right|\geq\frac{1}{2}\cdot\frac{1}{N}\sum_{k\in\left(\frac{N}{4},\frac{N}{2}\right]}\frac{1}{B}\left|c_{\left\lceil \frac{N}{2}\right\rceil }\right|
\]
\[
\geq\frac{1}{2N}\cdot\frac{N}{5}\cdot\frac{1}{B}\left|c_{\left\lceil \frac{N}{2}\right\rceil }\right|=\frac{1}{10B}\left|c_{\left\lceil \frac{N}{2}\right\rceil }\right|\geq\frac{1}{10B^{2}}\left|c_{n}\right|,
\]
in particular, $\left|c_{N}\right|\leq10B^{2}\left|a_{N}\right|$. 

Define 
\[
b_{n}=\begin{cases}
\begin{array}{cc}
a_{n} & \mbox{if }n\leq N\\
c_{n} & \mbox{if }n>N
\end{array} & ;\qquad d_{n}=\begin{cases}
\begin{array}{cc}
c_{n}-a_{n} & \mbox{if }n\leq N\\
0 & \mbox{if }n>N
\end{array} & .\end{cases}\end{cases}
\]
By Lemma 7.2,$\left\{ b_{n}\right\} \in GMS_{\alpha,\varphi}\left(3B+6B^{2}\cdot10B^{2}\right)\subset GMS_{\alpha,\varphi}\left(63B^{4}\right)$
and $b_{n}+d_{n}=c_{n}$, all $n$. Furthermore: 
\[
\left\Vert \left\{ b_{n}\right\} \right\Vert _{l_{\frac{1}{k}}^{1}}=\left(\sum_{n=1}^{N}+\sum_{n=N+1}^{\infty}\right)\frac{\left|b_{n}\right|}{n}
\]
\[
=\sum_{n=1}^{N}\frac{\left|a_{n}\right|}{n}+\sum_{n=N+1}^{\infty}\frac{\left|c_{n}\right|}{n}=\sigma_{N}+\sum_{n=N+1}^{\infty}\frac{\left|c_{n}\right|}{n};
\]
\[
\left\Vert \left\{ d_{n}\right\} \right\Vert _{l^{1}}=\sum_{n=1}^{N}\left|c_{n}-a_{n}\right|
\]
\[
\leq\sum_{n=1}^{N}\left|a_{n}\right|+\sum_{n=1}^{N}\left|c_{n}\right|=\frac{\sigma_{N}}{N}\sum_{n=1}^{N}n+\sum_{n=1}^{N}\left|c_{n}\right|
\]

\[
=\sigma_{N}\frac{N+1}{2}+N\sigma_{N}=\sigma_{N}\frac{3N+1}{2}.
\]

Therefore, observing that $t\leq{\displaystyle \frac{1}{N-1}}\leq{\displaystyle \frac{2}{N}}$,
we obtain: 
\[
\left\Vert \left\{ b_{n}\right\} \right\Vert _{l_{\frac{1}{k}}^{1}}+t\left\Vert \left\{ d_{n}\right\} \right\Vert _{l^{1}}
\]
\[
=\sigma_{N}\left(1+\frac{3N+1}{2}t\right)+\sum_{n=N+1}^{\infty}\frac{\left|c_{n}\right|}{n}\leq\sigma_{N}\left(1+\frac{3N+1}{N}\right)+\sum_{n=N+1}^{\infty}\frac{\left|c_{n}\right|}{n}
\]
\[
\leq\frac{9}{2}\sigma_{N}+\sum_{n=N+1}^{\infty}\frac{\left|c_{n}\right|}{n}=\frac{9}{2N}{\displaystyle \sum_{n=1}^{N}}\left|c_{n}\right|+\sum_{n=N+1}^{\infty}\frac{\left|c_{n}\right|}{n}
\]
\[
=\frac{9}{2N}\left({\displaystyle \sum_{n=1}^{\left\lfloor \frac{1}{t}\right\rfloor }}\left|c_{n}\right|+\left|c_{N}\right|\right)+\sum_{n=N+1}^{\infty}\frac{\left|c_{n}\right|}{n}
\]
\[
\leq\frac{9}{2}\left(t{\displaystyle \sum_{n=1}^{\left\lfloor \frac{1}{t}\right\rfloor }}\left|c_{n}\right|+\sum_{n=N}^{\infty}\frac{\left|c_{n}\right|}{n}\right)=\frac{9}{2}K\left(t,\left\{ c_{n}\right\} ,l_{\frac{1}{k}}^{1},l^{1}\right),
\]
in the last line we applied $N=\left\lfloor \frac{1}{t}\right\rfloor +1\in\left(\frac{1}{t},\frac{1}{t}+1\right]$$.\square$

B. Booton {[}3{]} gave a different proof of the following theorem
for $1<p$, $q\geq1$ and $\left\{ c_{k}\right\} \in GMS^{+}$.

$ $

\textbf{Theorem 7.4. }Let $1<p<\infty$, $0<q\leq\infty$. Assume
$\left\{ c_{k}\right\} \in GMS_{\alpha,\varphi}$. Let $f(x)={\displaystyle \sum_{k=1}^{\infty}}c_{k}e^{ikx}$.
Then $\exists C=C(p,q)$ such that
\begin{equation}
\left\Vert f\right\Vert _{L(p,q)(0,\pi)}\leq C\left\Vert \left\{ c_{n}\right\} \right\Vert _{l(p',q)}.\label{eq:Main result}
\end{equation}

\textbf{Proof.}

Let $X=l_{\frac{1}{k}}^{1}$ and let $w$ be the weight: $w(k)=k$.
Then $X_{w}=l^{1}$. Both $X$ and $X_{w}$ can be considered as quasi-normed
monoids, and so, by Gilbert's theorem (2.6), 
\[
\left\Vert \left\{ c_{k}\right\} \right\Vert _{\left(X,X_{w}\right){}_{\theta,q;K}}\sim\left(\intop_{0}^{\infty}\left(t^{-\theta}\left\Vert \left\{ c_{k}\right\} \cdot\sigma_{t}(k)\right\Vert _{X}\right)^{q}\frac{dt}{t}\right)^{1/q}
\]
\[
=\left(\intop_{0}^{\infty}\left(t^{-\theta}\sum_{k=1}^{\infty}tk\sigma(tk)\frac{\left|c_{k}\right|}{k}\right)^{q}\frac{dt}{t}\right)^{1/q}=\left(\intop_{0}^{\infty}\left(t^{1-\theta}\sum_{k=1}^{\infty}\sigma(tk)\left|c_{k}\right|\right)^{q}\frac{dt}{t}\right)^{1/q}
\]
\[
=\left(\intop_{0}^{\infty}\left(t^{\theta-1}\sum_{k=1}^{\infty}\sigma\left(\frac{k}{t}\right)\left|c_{k}\right|\right)^{q}\frac{dt}{t}\right)^{1/q}.
\]

Since $\sigma=I_{[1,2)}$ satisfies the conditions of Gilbert's theorem,
\[
\left(\intop_{0}^{\infty}\left(t^{\theta-1}\sum_{k=1}^{\infty}\sigma\left(\frac{k}{t}\right)\left|c_{k}\right|\right)^{q}\frac{dt}{t}\right)^{1/q}\sim\left(\intop_{0}^{\infty}\left(t^{\theta-1}\sum_{t\leq k<2t}\left|c_{k}\right|\right)^{q}\frac{dt}{t}\right)^{1/q}.
\]

$\left\{ c_{k}\right\} \in GMS_{1}\left(B\right)$ for some $B\geq1$,
and so: 
\[
{\displaystyle \sum_{t\leq k<2t}}\left|c_{k}\right|\leq\#\left\{ k:t\leq k<2t\right\} \cdot B\left|c_{\left\lceil t\right\rceil }\right|\leq2Bt\left|c_{\left\lceil t\right\rceil }\right|.
\]

Therefore, for $q<\infty$: 
\[
\left(\intop_{0}^{\infty}\left(t^{\theta-1}\sum_{t\leq k<2t}\left|c_{k}\right|\right)^{q}\frac{dt}{t}\right)^{1/q}\leq\left(\intop_{0}^{\infty}\left(2Bt^{\theta}\left|c_{\left\lceil t\right\rceil }\right|\right)^{q}\frac{dt}{t}\right)^{1/q}
\]
\[
=2B\left(\sum_{m=1}^{\infty}\intop_{m-1}^{m}\left(t^{\theta}\left|c_{\left\lceil t\right\rceil }\right|\right)^{q}\frac{dt}{t}\right)^{1/q}=2B\left(\sum_{m=1}^{\infty}\left|c_{m}\right|^{q}\intop_{m-1}^{m}t^{\theta q}\frac{dt}{t}\right)^{1/q}
\]
\[
\leq BC\left(\theta,q\right)\left(\sum_{m=1}^{\infty}\left|c_{m}\right|^{q}m^{\theta q-1}\right)^{1/q}=BC\left(\theta,q\right)\left\Vert \left\{ c_{k}\right\} \right\Vert _{l_{w\left(\frac{1}{\theta},q\right)}^{q}},
\]
where $C\left(\theta,q\right)=\max\left\{ 1,{\displaystyle \frac{1}{\theta q}},2^{1-\theta q}\right\} $.

Again, since $\left\{ c_{k}\right\} \in GMS_{1}\left(B\right)$, it
follows that for $t\leq k<2t$, $\left|c_{k}\right|\geq{\displaystyle \frac{1}{B}}\left|c_{\left\lfloor 2t\right\rfloor }\right|$.
Also:

If $t\leq{\displaystyle \frac{1}{2}}$ and $k\in\mathbb{N}$ then
$\sigma\left(\frac{k}{t}\right)=0$.

If $n-\frac{1}{2}<t\leq n$, where $n\in\mathbb{N}$, then 
\[
\#\left\{ k:t\leq k<2t\right\} =\#\left\{ k:n\leq k\leq2n-1\right\} =n>t.
\]

If $n<t\leq n+\frac{1}{2}$, where $n\in\mathbb{N}$, then 
\[
\#\left\{ k:t\leq k<2t\right\} =\#\left\{ k:n+1\leq k\leq2n\right\} =n>t-\frac{1}{2}>t-{\displaystyle \frac{t}{2}}=\frac{t}{2}.
\]

Therefore, for $t\geq\frac{1}{2}$, ${\displaystyle \sum_{t\leq k<2t}}\left|c_{k}\right|\geq{\displaystyle \frac{t}{2B}}\left|c_{\left\lfloor 2t\right\rfloor }\right|$,
and so: 
\[
\left(\intop_{0}^{\infty}\left(t^{\theta-1}\sum_{t\leq k<2t}\left|c_{k}\right|\right)^{q}\frac{dt}{t}\right)^{1/q}=\left(\intop_{\frac{1}{2}}^{\infty}\left(t^{\theta-1}\sum_{t\leq k<2t}\left|c_{k}\right|\right)^{q}\frac{dt}{t}\right)^{1/q}
\]
\[
\geq\left(\intop_{\frac{1}{2}}^{\infty}\left(\frac{t^{\theta}}{2B}\left|c_{\left\lfloor 2t\right\rfloor }\right|\right)^{q}\frac{dt}{t}\right)^{1/q}=\frac{1}{2B}\left(\sum_{m=1}^{\infty}\intop_{\frac{m}{2}}^{\frac{m+1}{2}}\left(t^{\theta}\left|c_{\left\lfloor 2t\right\rfloor }\right|\right)^{q}\frac{dt}{t}\right)^{1/q}
\]
\[
=\frac{1}{2B}\left(\sum_{m=1}^{\infty}\left|c_{m}\right|^{q}\intop_{\frac{m}{2}}^{\frac{m+1}{2}}t^{\theta q}\frac{dt}{t}\right)^{1/q}\geq\frac{C\left(\theta,q\right)}{B}\left(\sum_{m=1}^{\infty}\left|c_{m}\right|^{q}m^{\theta q-1}\right)^{1/q}
\]
\[
=\frac{C\left(\theta,q\right)}{B}\left\Vert \left\{ c_{k}\right\} \right\Vert _{l_{w\left(\frac{1}{\theta},q\right)}^{q}},
\]
where $C\left(\theta,q\right)=\min\left\{ \frac{1}{2}.\left(\frac{1}{2}\right)^{\theta q-1}\right\} $.

Therefore:
\begin{equation}
\left\Vert \left\{ c_{k}\right\} \right\Vert _{\left(X,X_{w}\right){}_{\theta,q;K}}\sim\left\Vert \left\{ c_{k}\right\} \right\Vert _{l_{w\left(\frac{1}{\theta},q\right)}^{q}}.\label{eq: last step in 8.4}
\end{equation}

On the other hand, Theorem 7.3 shows that 
\[
\left\Vert \left\{ c_{k}\right\} \right\Vert _{\left(X,X_{w}\right){}_{\theta,q;K}}=\left\Vert \left\{ c_{k}\right\} \right\Vert _{\left(l_{\frac{1}{k}}^{1},l^{1}\right)_{\theta,q;K}}\sim\left\Vert \left\{ c_{k}\right\} \right\Vert _{\left(l_{\frac{1}{k}}^{1}\cap GMS_{\alpha,\varphi},l^{1}\right)_{\theta,q;K}}.
\]

$l_{\frac{1}{k}}^{1}\cap GMS_{\alpha,\varphi}$ is a quasi-normed
monoid in the vector space of sequences (see Corollary 5.3), with
the norm $\left\Vert \cdot\right\Vert _{l_{\frac{1}{k}}^{1}}$. Assume
that $T$ is defined as in $(\ref{eq:Fourier series operator on GM - definition})$.
By Theorem 6.3, $T$ maps $G_{0}$ into $L(1,\infty)$, and $\left\Vert T\left\{ c_{n}\right\} \right\Vert _{L\left(1,\infty\right)}\leq6\pi B\left\Vert \left\{ c_{n}\right\} \right\Vert _{l_{\frac{1}{k}}^{1}}$.
Also, since a function with an absolutely convergent Fourier series
is bounded, it follows that $T:G_{1}\rightarrow L(\infty,\infty)$.
Therefore, 
\[
\left\Vert T\left\{ c_{n}\right\} \right\Vert _{\left(L\left(1,\infty\right),L\left(\infty,\infty\right)\right)\theta,q;K}\leq C\left(B\right)\left\Vert \left\{ c_{n}\right\} \right\Vert _{\left(G_{0},G_{1}\right){}_{\theta,q;K}}\sim\left\Vert \left\{ c_{n}\right\} \right\Vert _{\left(X,X_{w}\right){}_{\theta,q;K}}.
\]

Finally, by the interpolation theorem,$\forall\theta\in(0,1),\forall q\in(0,\infty)$:
\[
(L(1,\infty),L(\infty,\infty))_{\theta,q}=L(p,q),
\]
where $p={\displaystyle \frac{1}{1-\theta}}$. That is to say, $\theta={\displaystyle \frac{1-p}{p}}={\displaystyle \frac{1}{p'}}$.

Since $\left\{ c_{k}\right\} \in GMS$ the conditions of Lemma 4.4
hold and, therefore, $\left|c_{n}\right|\leq Bc_{\left\lfloor \frac{n}{2}\right\rfloor }^{*}$.
This implies that $\left\Vert \left\{ c_{k}\right\} \right\Vert _{l_{w(p',q)}^{q}}\leq B'(p,q)\left\Vert \left\{ c_{k}^{*}\right\} \right\Vert _{l_{w(p',q)}^{q}}=B'(p,q)\left\Vert \left\{ c_{k}\right\} \right\Vert _{l(p,q)}$,
thus obtaining $(\ref{eq:Main result})$.$\square$

Combining the results of Theorems 7.4 and 6.4, we obtain the main
result of this section.

$ $

\textbf{Theorem 7.5}. Let $1<p<\infty$, $0<q\leq\infty$. Assume
$\left\{ c_{k}\right\} \in GMS_{\alpha,\varphi}$. Define $f(x)={\displaystyle \sum_{k=1}^{\infty}}c_{k}e^{ikx}$.
Then: 
\begin{equation}
\left\Vert \left\{ c_{n}\right\} \right\Vert _{l(p',q)}\sim\left\Vert f\right\Vert _{L(p,q)(0,\pi)}.\label{eq:Main result-1-1}
\end{equation}

\textbf{8. On Hardy's inequality for general monotone functions.}

$ $

\textbf{Lemma 8.1}. If $f,g\in GM$ then $f\cdot g\in GM$.

\textbf{Proof}.

\[
V_{fg}\left([x,2x]\right)=\sup_{\left\{ x_{j}\right\} \in{\cal P}([x,2x])}\left\{ \sum_{j=0}^{n-1}\left|f\left(x_{j+1}\right)g\left(x_{j+1}\right)-f\left(x_{j}\right)g\left(x_{j}\right)\right|\right\} 
\]
\[
\leq\sup_{\left\{ x_{j}\right\} \in{\cal P}([x,2x])}\left\{ \sum_{j=0}^{n-1}\left|f\left(x_{j+1}\right)\right|\left|g\left(x_{j+1}\right)-g\left(x_{j}\right)\right|\right\} 
\]
\[
+\sup_{\left\{ x_{j}\right\} \in{\cal P}([x,2x])}\left\{ \sum_{j=0}^{n-1}\left|g\left(x_{j}\right)\right|\left|f\left(x_{j+1}\right)-f\left(x_{j}\right)\right|\right\} 
\]
\[
\leq\sup_{[x,2x]}\left|f(t)\right|\cdot V_{g}\left([x,2x]\right)+\sup_{[x,2x]}\left|g(t)\right|\cdot V_{f}\left([x,2x]\right).
\]

Since $f,g\in GM$ then $\exists B_{1}>0$, $B_{2}>0$ such that $\forall x>0$:
\[
\begin{cases}
\begin{array}{c}
V_{f}\left([x,2x]\right)\leq B_{1}\left|f(x)\right|\\
V_{g}\left([x,2x]\right)\leq B_{2}\left|g(x)\right|
\end{array} & .\end{cases}
\]
Also, by $(\ref{eq:almost monotone-1: B}),$ for all $x>0$ and $t\in[x,2x]$:
\[
\begin{cases}
\begin{array}{c}
\left|f(t)\right|\leq2B_{1}\left|f(x)\right|\\
\left|g(t)\right|\leq2B_{2}\left|g(x)\right|
\end{array} & .\end{cases}
\]
Therefore, $\forall x>0$: 
\[
V_{fg}\left([x,2x]\right)\leq4B_{1}B_{2}\left|f(x)g(x)\right|,
\]
that is, $fg\in GM(4B_{1}B_{2}).\square$

$ $

We use interpolation to give a very short proof of Hardy's inequality
for general monotone functions, for $\forall\alpha>0$, $\forall q>0$.

$ $

$ $

\textbf{Theorem 8.2} (Hardy's inequality for general monotone functions).

Let $\alpha>0$, $0<q\leq\infty$, and $f\in GM^{+}$. Then 
\begin{equation}
\left(\intop_{0}^{\infty}\left(x^{-\alpha}\intop_{0}^{x}f(t)\frac{dt}{t}\right)^{q}\frac{dx}{x}\right)^{\frac{1}{q}}\leq C(\alpha,q)\left(\intop_{0}^{\infty}\left(x^{-\alpha}f(x)\right)^{q}\frac{dx}{x}\right)^{\frac{1}{q}}.\label{eq:Hardy for GM}
\end{equation}

\textbf{Proof.}

1) Consider $\alpha\leq1$. Take $g(t)={\displaystyle \frac{f(t)}{t}}$.
By Lemma 8.1, since ${\displaystyle \frac{1}{t}}\in GM^{+}$, $g\in GM^{+}$.

Let $T$ be defined for integrable on $[0,a]$ functions, for all
$a>0$, as 
\[
\left(Th\right)(x)=\frac{1}{x}\intop_{0}^{x}h(t)dt.
\]

Then $T:L^{\infty}\rightarrow L^{\infty}$ and $T:L(1,1)\rightarrow L(1,\infty)$.
Applying the interpolation theorem, we obtain

\[
T:L(p,q)\rightarrow L(p,q),
\]
where $1<p<\infty$, $0<q\leq\infty$.

This implies 
\[
\left\Vert Tg\right\Vert _{L(p,q)}=\left(\intop_{0}^{\infty}x^{\frac{q}{p}}\left(\left(\frac{1}{x}\intop_{0}^{x}g(t)dt\right)^{*}\right)^{q}\frac{dx}{x}\right)^{\frac{1}{q}}
\]
\begin{equation}
\leq C(p,q)\left(\intop_{0}^{\infty}x^{\frac{q}{p}}\left(g^{*}(x)\right)^{q}\frac{dx}{x}\right)^{\frac{1}{q}}=C(p,q)\left\Vert g\right\Vert _{L(p,q)}.\label{eq:Interpolation of the mean}
\end{equation}

By Theorem 4.10, $\left\Vert g\right\Vert _{L(p,q)}\sim\left\Vert g\right\Vert _{L_{w(p,q)}^{q}}$.
By Lemma 5.6, $Tg\in GM^{+}$, hence $\left\Vert Tg\right\Vert _{L(p,q)}\sim\left\Vert Tg\right\Vert _{L_{w(p,q)}^{q}}$.
Therefore, $(\ref{eq:Interpolation of the mean})$ can be written
as 
\[
\left(\intop_{0}^{\infty}x^{\frac{q}{p}}\left(\frac{1}{x}\intop_{0}^{x}g(t)dt\right)^{q}\frac{dx}{x}\right)^{\frac{1}{q}}\sim\left\Vert Tg\right\Vert _{L(p,q)}
\]
\[
\leq C(p,q)\left\Vert g\right\Vert _{L(p,q)}\sim C(p,q)\left(\intop_{0}^{\infty}x^{\frac{q}{p}}\left(g(x)\right)^{q}\frac{dx}{x}\right)^{\frac{1}{q}}.
\]

Since $g(t)={\displaystyle \frac{f(t)}{t}}$, we obtain

\begin{equation}
\left(\intop_{0}^{\infty}x^{\frac{q}{p}-q}\left(\intop_{0}^{x}f(t)\frac{dt}{t}\right)^{q}\frac{dx}{x}\right)^{\frac{1}{q}}\leq C(p,q)\left(\intop_{0}^{\infty}x^{\frac{q}{p}-q}\left(f(x)\right)^{q}\frac{dx}{x}\right)^{\frac{1}{q}}.\label{eq:Hardy in the interpolation form}
\end{equation}

Inequality $(\ref{eq:Hardy for GM})$ for $\alpha\leq1$ follows from
$(\ref{eq:Hardy in the interpolation form})$ if we take $\alpha=1-\frac{1}{p}=\frac{1}{p'}$.

2) Consider $\alpha>1$. Take $g(t)=t^{\epsilon-\alpha}f(t)$, where
$0<\epsilon<1$. Clearly, $g\in GM$. Also: 
\[
x^{-\alpha}\intop_{0}^{x}f(t)\frac{dt}{t}=x^{-\epsilon}\cdot\intop_{0}^{x}x^{\epsilon-\alpha}f(t)\frac{dt}{t}\leq x^{-\epsilon}\cdot\intop_{0}^{x}g(t)\frac{dt}{t}.
\]
By $(\ref{eq:Hardy for GM})$, for $\epsilon<1$ and $g\in GM$, 
\[
\left(\intop_{0}^{\infty}\left(x^{-\epsilon}\intop_{0}^{x}g(t)\frac{dt}{t}\right)^{q}\frac{dx}{x}\right)^{\frac{1}{q}}\leq C(\alpha,q)\left(\intop_{0}^{\infty}\left(x^{-\epsilon}g(x)\right)^{q}\frac{dx}{x}\right)^{\frac{1}{q}}
\]
\[
=C(\alpha,q)\left(\intop_{0}^{\infty}\left(x^{-\alpha}f(x)\right)^{q}\frac{dx}{x}\right)^{\frac{1}{q}},
\]
and $(\ref{eq:Hardy for GM})$ for $\alpha>1$ follows.$\square$

\[
\]

\textbf{References}.

{[}1{]} Askey, R.; Wainger, S. Integrability theorems for Fourier
series. Duke Math. J. 33, 223\textendash 228, 1966.

{[}2{]} Bergh, J.; Löfström, J. Interpolation spaces. An introduction.
Grundlehren der Mathematischen Wissenschaften, No. 223. Springer-Verlag,
Berlin-New York, 1976. x+207 pp.

{[}3{]} Booton, B. General monotone sequences and trigonometric series.
Math. Nachr. 1-12, 2013.

{[}4{]} Brudny\u{\i}, Yu. A.; Krugljak, N. Ya. Interpolation functors
and interpolation spaces. Vol. I. North-Holland Mathematical Library,
47. North-Holland Publishing Co., Amsterdam, 1991. xvi+718 pp.

{[}5{]} Gilbert, J. E. Interpolation between weighted $L^{p}$-spaces.
Ark. Mat. 10, 235\textendash 249, 1972.

{[}6{]} Hardy, G.H.; Littlewood, J.E. Some new properties of Fourier
constants. J.L.M.S. 6, 3-9, 1931.

{[}7{]} Hardy, G.H.; Littlewood, J.E. Some new properties of Fourier
constants. M.A. 97, 159-209, 1926.

{[}8{]} Hardy, G.H.; Littlewood, J.E.; Pólya, G. Inequalities. 2d
ed. Cambridge, at the University Press, 1952.

{[}9{]} Liflyand, E.; Tikhonov, S. A concept of general monotonicity
and applications. Math. Nachr. 284, no. 8-9, 1083\textendash 1098,
2011.

{[}10{]} Paley, R.E.A.C. \textquotedbl{}Some theorems on orthogonal
functions (1).\textquotedbl{} Studia Mathematica 3.1: 226-238, 1931.

{[}11{]} Sagher, Y. An application of interpolation theory to Fourier
series. Studia Math. 41, 169\textendash 181, 1972.

{[}12{]} Sagher, Y. Integrability conditions for the Fourier transform.
J. Math. Anal. Appl. 54, no. 1, 151\textendash 156, 1976.

{[}13{]} Shah, S. M. A note on quasi-monotone series. Math. Student
15 (1947), 19\textendash 24, 1948.

{[}14{]} Tikhonov, S. Trigonometric series with general monotone coefficients.
J. Math. Anal. Appl. 326, no. 1, 721\textendash 735, 2007.

{[}15{]} Zygmund, A. Trigonometric series. 2nd ed. Vols. I, II. Cambridge
University Press, New York 1959. Vol. I. xii+383 pp.; Vol. II. vii+354
pp.
\end{document}